\documentclass[MSNbibl,number,citesort,seceqn,dvips]{arxbj}
\usepackage{upgreek,dcolumn,url,breakurl}
\usepackage[left]{vmkcol}

\aid{0}
\volume{20}
\issue{4}
\pubyear{2014}
\firstpage{1738}
\lastpage{1764}
\doi{10.3150/13-BEJ539} %kopijuoti is PTS

\makeatletter
\newcolumntype{d}[1]{D{.}{.}{#1}}

\newcommand{\rrvert}{\vert}
\newcommand{\llvert}{\vert}
\newtheorem{theorem}{Theorem}[section]

\newtheorem{proposition}{Proposition}[section]
\newremark{example}{Example}%[section]
\newtheorem{lemma}[theorem]{Lemma}
\newremark{remark}{Remark}[section]

\def\tr{\operatorname{tr}}
\def\l{\ell}
\def\p{\partial}
\def\zero{\mathbf{0}}
\def\one{\mathbf{1}}
\def\diag{\operatorname{diag}}
\def\E{\mathrm{E}}
\def\Cov{\operatorname{Cov}}
\def\rt#1{\sqrt{#1}}
\def\oo{\infty}
\def\vec{\operatorname{vec}}
\def\vech{\operatorname{vech}}
\def\bksl{$\setminus$}

\def\Amat{\mathbf{A}}
\def\bbf{\mathbf{b}}
\def\Bmat{\mathbf{B}}
\def\Cmat{\mathbf{C}}
\def\Dmat{\mathbf{D}}
\def\Emat{\mathbf{E}}
\def\H{\mathbf{H}}
\def\Imat{\mathbf{I}}
\def\J{\mathbf{J}}
\def\Kmat{\mathbf{K}}
\def\Smat{\mathbf{S}}
\def\Tmat{\mathbf{T}}
\def\V{\mathbf{V}}
\def\xbf{\mathbf{x}}
\def\Ybf{\mathbf{Y}}
\def\Zbf{\mathbf{Z}}

\def\b{\beta}
\def\bebf{\bolds{\beta}}
\def\gbf{\bolds{\gamma}}
\def\ghatbf{\hat{\gbf}}
\def\dbf{\bolds{\delta}}
\def\epsbf{\bolds{\varepsilon}}
\def\th{\theta}
\def\thbf{\bolds{\theta}}
\def\thatbf{\hat{\thbf}}
\def\thtilbf{\tilde{\thbf}}
\def\vth{\vartheta}
\def\vthbf{\bolds{\vartheta}}
\def\ze{\zeta}
\def\zebf{\bolds{\zeta}}
\def\zehatbf{\hat{\zebf}}
\def\lmbf{\bolds{\lambda}}
\makeatother

\begin{document}
\begin{frontmatter}

\title{Model comparison with composite likelihood information criteria}

\runtitle{Model comparison with composite likelihood}

\begin{aug}
%%%% inicialai - be tarpu
\author[a,b]{\inits{C.T.}\fnms{Chi Tim} \snm{Ng}\thanksref{a,b}\ead[label=e1,mark]{easterlyng@gmail.com}}
\and
\author[c]{\inits{H.}\fnms{Harry} \snm{Joe}\thanksref{c}\ead[label=e2]{Harry.Joe@ubc.ca}}
%%\runauthor{} %% auto
\address[a]{Department of Statistics, Seoul
National University, Room 430, Building 25, Seoul,\\ South Korea.
\printead{e1}}

\address[b]{Department of Statistics, Chonnam National University,
Gwangju, 500-757, South Korea}

\address[c]{Department of Statistics, University of British Columbia, Room ESB 3138, Earth Sciences
Building, Vancouver, Canada.
\printead{e2}}
\end{aug}

% HISTORY:
\received{\smonth{11} \syear{2012}}
\revised{\smonth{6} \syear{2013}}

% ABSTRACT
%
\begin{abstract}
Comparisons are made for the amount of agreement of the composite
likelihood information criteria and their full likelihood counterparts when
making decisions among the fits of different models, and some properties
of penalty term for composite likelihood information criteria are
obtained. Asymptotic theory is given for
the case when a simpler model is nested within a bigger model, and
the bigger
model approaches the simpler model under a sequence of local alternatives.
Composite likelihood can more or less frequently choose the bigger model,
depending on the direction of local alternatives; in the former case,
composite likelihood has more ``power'' to choose the bigger model.
The behaviors of the
information criteria are illustrated via theory and simulation examples
of the Gaussian linear mixed-effects model.
\end{abstract}

% KEYWORDS
% visi is mazosios raides ir pagal abecele
%

\end{frontmatter}

%s1 #&#
\section{Introduction}

Composite likelihood inference based on low-dimensional marginal or
conditional distributions is common when the full likelihood is
computationally too difficult. It has been increasing used in recent
years for inference with complex models; see Varin \cite{Va}, Varin
\textit{et al.} \cite{VRF} for reviews.

For model selection with composite likelihood, one might wonder if the
use of limited or reduced information leads to different decisions. To
understand this,
an asymptotic theory based on the theory of a sequence of contiguous
local alternatives is developed to compare Akaike information criterion
(AIC) and Bayesian information criterion (BIC) in their full likelihood
and composite marginal likelihood versions.
We show that model selection based on AIC and its composite likelihood
counterpart (as proposed in Varin and Vidoni \cite{VV}) are sometimes
similar (when
models under consideration are far apart) and sometimes not similar
(when one model is a perturbation of another). The patterns can be
explained via local alternatives where the perturbed model is at a
distance $n^{-1/2}$ from a ``null'' or simplified model, with $n$ being
the sample size.

We also provide simulation results under models where the maximum
likelihood is feasible; one class of such models is the linear
mixed-effects models based on the normal distribution. Within different
sub-cases of the Gaussian linear mixed-effects models, the simulation
results are consistent with the asymptotic theory.

The remainder of the paper is organized as follows. In Section~\ref{sec-clic}, we introduce our notation and state the definitions for the
composite marginal likelihood and the information criteria. In Section~\ref{sec-theorems}, asymptotic properties of composite likelihood
information criteria are presented. In Section~\ref{sec-simulation},
comparisons of decisions between Varin and Vidoni's composite
likelihood information criterion (abbreviated CLAIC as
in Varin \textit{et al.} \cite{VRF}), Gao and Song's information criterion
(abbreviated as
CLBIC in Gao and Song \cite{GS}), and their full-likelihood
counterparts are
summarized via simulation studies.
Section~\ref{sec-realdata} contains a data example with a
mixed-effects model.
Section~\ref{sec-discussion} concludes with some discussion and future
research. The proofs of the main theorems in Section~\ref{sec-theorems} are given in
Appendix \ref{appa}.

%s2 #&#
\section{Composite likelihood and information criteria}
\label{sec-clic}

For the comparison of composite likelihood and full likelihood
information criteria, we consider the case of independent multivariate
measurements on $n$ subjects, possibly
with covariates. Nested statistical models will be considered.

%s2.1 #&#
\subsection{Model}

Let $\mathbf{y}_1,\ldots,\mathbf{y}_n$ be the realizations of independent
$d$-dimensional random vectors $\Ybf_i$, with respective covariates summarized
as matrices $\xbf_1,\ldots,\xbf_n$. Suppose that the data generating
mechanism of $\Ybf_i$ is governed by the density function $g(\mathbf{
y}_i;\xbf
_i)$. Candidate parametric
models are $f^{(M)}(\mathbf{y}_i;\xbf_i,\thbf^{(M)})$, for
$M=1,2,\ldots,$;
$M$ is an index for different models that are considered, and $\thbf
^{(M)}$ is the parameter vector for
model $M$. Let $p_M=\dim(\thbf^{(M)})$ be the dimension of $\thbf
^{(M)}$ for a generic model $M$; the superscript will be omitted unless
we are referring to two or more models.

%s2.2 #&#
\subsection{Composite likelihood}

For model $M$, let $L_{\mathrm{CL}}^{(M)}(\thbf^{(M)})=
L_{\mathrm{CL}}^{(M)}(\thbf^{(M)};\mathbf{y}_1,\ldots,\mathbf{y}_n;\xbf
_1,\ldots
,\xbf_n)$ be
a particular composite
marginal log-likelihood. We are using the same composite likelihood
(same set of marginal density functions) for all competing models.
Let $S\subset\{1,2,\ldots,d\}$ be a non-empty subset of indexes.
For notation, $f_S^{(M)}$ indicates a marginal density of $f^{(M)}$
with margin
$S$ and $g_S$ is the corresponding margin of $g$.
The particular composite likelihood could be based on all bivariate
margins, or a subset of bivariate margins,
or more generally a set of margins $\{S_1,\ldots,S_Q\}$
with corresponding weights $w_1,\ldots,w_Q$.
Suppressing the superscript for the model, let
%
%e2.1 #&#
%
\begin{equation}
L_{\mathrm{CL}}(\thbf)=L_{\mathrm{CL}}(\thbf;\mathbf{y}_1,\ldots,
\mathbf{y}_n;\xbf _1,\ldots,\xbf_n)= \sum
_{i=1}^n \l_{\mathrm{CL}}(\thbf;
\mathbf{y}_i,\xbf_i) \label{eq-CL}
\end{equation}
be the log composite likelihood. Here
%
%e2.2 #&#
%
\begin{equation}
\exp\bigl\{\l_{\mathrm{CL}}(\thbf;\mathbf{y}_i,
\xbf_i)\bigr\}= \prod_{q=1}^Q
f_{S_q}^{w_q}(\mathbf{y}_{i,S_q};\xbf_i,
\thbf), \label{eq-cml}
\end{equation}
$S_q$ is a subset consisting of indexes, and $w_q$ is a positive weight
for $S_q$.
For example, if these are the pairs for bivariate composite likelihood,
then the cardinality of $\{S_q\}$ is $Q=d(d-1)/2$.
Note that the case of full likelihood is covered with
$S_1=\{1,\ldots,d\}$ with the cardinality of $\{S_q\}$ being~1.

%s2.3 #&#
\subsection{Composite likelihood information criteria}

Consider the composite likelihood versions of
Akaike information criterion (AIC) and Bayesian information criterion
(BIC) described in Varin and Vidoni \cite{VV}, Gao and Song \cite
{GS}, Varin \textit{et al.}~\cite{VRF}. They are defined as (with
superscript for model $M$ omitted):
%
%e2.3 #&#
%
\begin{equation}
\mathrm{\mathrm{CLAIC}}=-2 L_{\mathrm{CL}}(\hat{\bolds{\theta}}_{\mathrm{CL}}) + 2 {
\rm tr}\bigl\{\mathbf{J}(\hat{\bolds{\theta}}_{\mathrm{CL}})
\mathbf{H}^{-1}(\hat{\bolds{\theta}}_{\mathrm{CL}})\bigr\}
\label{eq-cl-aic}
\end{equation}
and
%
%e2.4 #&#
%
\begin{equation}
\mathrm{\mathrm{CLBIC}}=-2 L_{\mathrm{CL}}(\hat{\bolds{\theta}}_{\mathrm{CL}}) + (
\log n) \operatorname{tr}\bigl\{\mathbf{J}(\hat{\bolds{\theta}}_{\mathrm{CL}})
\mathbf{H}^{-1}(\hat{\bolds{\theta}}_{\mathrm{CL}})\bigr\}.
\label{eq-cl-bic}
\end{equation}
Here, $\thatbf_{\mathrm{CL}}=\thatbf_{n,\mathrm{CL}}$ is the composite likelihood estimator
that maximizes (\ref{eq-CL}).
The matrices $\H(\thbf)$ and $\J(\thbf)$ are the Hessian matrix and the
covariance matrix of the score function, respectively,
\[
\H(\thbf)=-\lim_{n\to\oo} n^{-1} {\p^2 L_{\mathrm{CL}}(\thbf;\mathbf{y}_1,\ldots,\mathbf{y}_n,\xbf_1,\ldots
,\xbf_n)
\over\p\thbf\,\p\thbf^T}
\]
and
\[
\J(\thbf)=\Cov \biggl[ n^{-1/2} {\p L_{\mathrm{CL}}(\thbf;\mathbf{y}_1,\ldots,\mathbf{y}_n,\xbf_1,\ldots
,\xbf_n)
\over\p\thbf} \biggr].
\]

When there are several models, the CLAIC (CLBIC) principle selects the
model with smallest
value of CLAIC (CLBIC). CLAIC has penalty term $2\tr(\J\H^{-1})$ and
CLBIC has penalty term $(\log n) \tr(\J\H^{-1})$ that depends on the sample
size $n$. With large $n$, CLBIC might choose smaller models than CLAIC.

%s3 #&#
\section{Main theorems}
\label{sec-theorems}

The main results are presented in this section, with proofs in the
\hyperref[app]{Appendix}. Consider the nested cases where model 1 is nested within
model 2. Proposition \ref{lik-ratio} gives general results of the
composite likelihood ratio under nested cases. If the true model is
covered by either model 1 or model 2, Theorem \ref{main} provides
further comparison of the asymptotic properties of CLAIC and CLBIC
under a sequence of local alternative hypotheses. Results under model
misspecification are summarized in Theorem \ref{main.robust}.

To describe the theorems, the following notation is used,
\begin{itemize}
\item
Model 1: $\Ybf|\xbf\sim f^{(1)}(\mathbf{y};\xbf,\thbf)$, $\thbf
\in
\Theta$.
\item
Model 2: $\Ybf|\xbf\sim f^{(2)}(\mathbf{y};\xbf,\gbf)$, $\gbf\in
\Gamma$.
\item
True model: $\Ybf|\xbf\sim g(\mathbf{y};\xbf)$.
\end{itemize}
This notation matches $\thbf^{(1)}=\thbf$ and $\thbf^{(2)}=\gbf$,
as used in Section~\ref{sec-clic}, but we are temporarily reducing the
number of superscripts. Let $\thbf^*$ be the parameters for
$f^{(1)}(\cdot;\thbf)$ such that $f^{(1)}$ is the closest to $g$ in the
divergence (see Xu and Reid \cite{XR})
based on the composite log-likelihood function $L^{(1)}_{\mathrm{CL}}$.
Similarly $\gbf^*$ is defined. Note that $\thbf^*$ and $\gbf^*$ might
depend on the composite log-likelihood that is used.

%pr3.1 #&#
\begin{proposition}[(Asymptotic distribution of the composite likelihood ratio)]\label{lik-ratio}
 Consider
the log composite likelihood ratio of two competing models,
%
%e3.1 #&#
%
\begin{equation}
\operatorname{LR}=L^{(2)}_{\mathrm{CL}}(\ghatbf)-L^{(1)}_{\mathrm{CL}}(
\thatbf). \label{eq-compLR}
\end{equation}
Suppose that assumptions \textup{A1--A3} (given in Appendix \textup{\ref{appa1}}) hold. If for
all $(\xbf,\mathbf{y}) $,
%
%e3.2 #&#
%
\begin{equation}
f^{(1)}\bigl(\mathbf{y};\xbf,\thbf^*\bigr) =f^{(2)}\bigl(
\mathbf{y};\xbf,\gbf^*\bigr), \label{eq-compden}
\end{equation}
then the limiting distribution of $2 \operatorname{LR}$ has the same law as
$\Zbf^TD \Zbf$, where $\Zbf$ is a vector of independent standard
normal random variables and $D$ is a diagonal matrix with eigenvalues
of the matrix:
%
%e3.3 #&#
%
\begin{equation}
\Bmat= \pmatrix{
-\bigl(\J^{(11)}\bigr)
\bigl(\H^{(1)}\bigr)^{-1}&\bigl(\J^{(12)}\bigr) \bigl(
\H^{(2)}\bigr)^{-1}
\vspace*{2pt}\cr
-\bigl(\J^{(21)}\bigr) \bigl(\H^{(1)}\bigr)^{-1}&
\bigl(\J^{(22)}\bigr) \bigl(\H^{(2)}\bigr)^{-1}}. \label{eq-Bmat}
\end{equation}
Here, $\H^{(1)}$, $\J^{(12)}$, etc., are defined in Appendix \textup{\ref{appa1}}.
\end{proposition}

In order to understand how different criteria can differ, we do an analysis
for a sequence of contiguous alternatives, in which the true model is
model 2
and its parameter depends on the
sample size $n$ and is closer to the null model as $n$ increases.
Such theory helps to explain what happens in finite samples;
see Section~\ref{sec-simulation}. Suppose that model 2 is
$f^{(2)}(\cdot
;\thbf,\zebf)$ and model~1 (null model) is nested within model 2, that
is, $f^{(1)}(\cdot;\thbf)=f^{(2)}(\cdot;\thbf,\zero)$.
The local alternatives assumption refers to that $g(\cdot)=f^{(2)}(\cdot
;\thbf^*_{2n},\zebf^*_n)$ with $\zebf^*_n=a_n\epsbf$ converges to
$\zebf^*=\zero$ at rate $a_n=n^{-1/2}$ or $a_n=\sqrt{\log n/n}$, and
$\thbf^*_{2n}\to\thbf^*$.
Let $\thbf^*_{1n}$ be the parameter
for $f^{(1)}(\cdot;\thbf)$ such\vadjust{\goodbreak} that $f^{(1)}$ is closest to $g$ in the
divergence
(see Xu and Reid \cite{XR}) based on the composite log-likelihood
function $L^{(1)}_{\mathrm{CL}}$.
Assume that $\thbf^*_{1n}$ and $\thbf^*_{2n}$ are asymptotically
equivalent, that is,
%
%e3.4 #&#
%
\begin{equation}
\thbf^*_{1n}-\thbf^*\to\zero,\qquad n\to\oo. \label{eq-locallimit}
\end{equation}
We next state the main theorem for comparing CLAIC, AIC, CLBIC
for nested models, when the null model is true, or when the larger model
is true under a sequence of local alternatives.

%th3.1 #&#
\begin{theorem}\label{main}
Consider the model selection problem $H_1$: Model 1 $f^{(1)}(\cdot
;\thbf
)$ is the true model versus $H_2$: Model 2 $f^{(2)}(\cdot;\thbf,\zebf)$
is the true model. Here, $\gbf$ is $p_2$-dimensional and $\zebf$ is
$m$-dimensional, where $m=p_2-p_1$. Let $P^{\mathrm{AIC}}_1$ be the probability
that AIC selects model 1. Similar notation is used for BIC, CLAIC, and CLBIC.

\begin{longlist}[(1)]
\item[(1)] Under $H_1$, $P^{\mathrm{CLAIC}}_1\to C_1\in(0,1)$ and
$P^{\mathrm{CLBIC}}_1\to1$.\vadjust{\goodbreak}

\item[(2)] Under $H_1$, $P^{\mathrm{CLAIC}}_1<P^{\mathrm{AIC}}_1$.

\item[(3)] Under $H_2$ with $\zebf=\zebf^*_n=\epsbf n^{-1/2}$ and
$\epsbf=\mathrm{O}(1)$, and assuming (\ref{eq-locallimit}), $P^{\mathrm{CLAIC}}_1\to
C_2\in(0,1)$ and $P^{\mathrm{CLBIC}}_1\to1$.

\item[(4)] Under $H_2$ with $\zebf=\zebf^*_n=\epsbf\sqrt{\log n/n}$
and $\epsbf=\mathrm{O}(1)$, and assuming (\ref{eq-locallimit}),
$P^{\mathrm{CLAIC}}_1\to0$ and $P^{\mathrm{CLBIC}}_1\to C_3\in(0,1)$.
\end{longlist}

 To be more specific, we have $C_1=P(\lambda_1 U_1+\cdots
+\lambda_m U_m<2(\lambda_1+\cdots+\lambda_m))$, where
$U_1,U_2,\ldots
,U_m$ are independent $\chi^2_1$ random variables and $\lambda
_1,\lambda
_2,\ldots,\lambda_m$ are the non-zero eigenvalues of $\Bmat$ defined in
(\ref{eq-Bmat}). If the full-likelihood is used, $\lambda_1=\cdots
=\lambda_m=1$.
\end{theorem}

In Theorem \ref{main}, (1) is a special case of (3) with $\epsbf
=\zero$.
The asymptotic results (1) and (3) are natural. Intuitively speaking,
if less parameters than the true model are selected, the composite
likelihood decreases by a positive quantity of $\mathrm{O}(n)$. Such a
decrease dominates the CLAIC (CLBIC) penalty term so the penalty term
is ignorable. This guarantees that the true model is better than the
smaller models in terms of CLAIC (CLBIC). On the other hand, if more
parameters are involved than necessary, the increase in composite
likelihood is just $\mathrm{O}(1)$. For CLAIC, the change in penalty term is
also $\mathrm{O}(1)$, so the model is correctly selected only with some
positive probability. For CLBIC, provided that the penalty term is
monotonic (see Lemma \ref{mono}), it is guaranteed that the change in
penalty term is positive and is $\mathrm{O}(\log n)$, dominating the increase
in composite likelihood. Then, the true model is better than any other
bigger model.

If model 2 is the true model and the two models are sufficiently far
apart from each other, that is, $\zebf=\mathrm{O}(1)\neq\zero$, then all the
criteria asymptotically choose the correct model. On the contrary, if
the two models differ by only a small perturbation, for example, $\zebf
=\mathrm{O}(1/\sqrt{n})$ or $\zebf=\mathrm{O}(\sqrt{\log n/n})$, it can be seen from
results (3) and (4) that the behavior of CLAIC and CLBIC differ. CLBIC
is less likely select the correct model than CLAIC.

Comparing CLAIC and its full-likelihood counterpart, CLAIC has greater
probability of selecting the larger model. The difference in such
probabilities depends on the eigenvalues $\lambda_1,\ldots,\lambda
_m
$. Roughly speaking, if $(\lambda_1,\ldots,\lambda_m)$, after
standardization is closer to $(1,\ldots,1)$, the ``loss of
information''\vadjust{\goodbreak} due to the use of composite likelihood is less
significant. It is natural to consider $C_1$ in Theorem \ref{main} as a
measurement of closeness of the composite likelihood to the
full-likelihood. It is interesting to note that $C_1$ does not depend
on the parameters for full-likelihood. For composite likelihood, it is
possible that $C_1$ depends on the parameters through $\lambda
_1,\ldots
,\lambda_m$. The dependence of $C_1$ on the parameters will be
illustrated via simulation examples in Section~\ref{sec-simulation}.

Part of the results in Theorem \ref{main} can be generalized to the
situation of model misspecification.

%th3.2 #&#
\begin{theorem}[(The same notation as in Proposition \ref{lik-ratio}
and Theorem \ref{main} is used)]\label{main.robust} Suppose that model 1 is nested within
model 2 but neither model 1 nor model 2 is the true model.
Let $(\thbf^*_{2n},\zebf^*_n)$ be the parameter under model 2 that is
the closest to the true model in the divergence (see Xu and Reid \cite
{XR}) based
on the composite likelihood.

 If equation (\ref{eq-compden}) holds, (1) $P^{\mathrm{CLAIC}}_1\to
C_1\in(0,1)$ and $P^{\mathrm{CLBIC}}_1\to1$.

 If equation (\ref{eq-compden}) does not hold,

(2) If $\zebf^*_n=\mathrm{O}(1)$,
then $P^{\mathrm{CLAIC}}_1\to0$ and $P^{\mathrm{CLBIC}}_1\to0$.

(3) If $\zebf^*_n=\epsbf n^{-1/2}$ and $\epsbf=\mathrm{O}(1)$, and
assuming (\ref{eq-locallimit}), $P^{\mathrm{CLAIC}}_1\to C_2\in(0,1)$ and
$P^{\mathrm{CLBIC}}_1\to1$.\vadjust{\goodbreak}

(4) If $\zebf^*_n=\epsbf\sqrt{\log n/n}$ and $\epsbf
=\mathrm{O}(1)
$, and assuming (\ref{eq-locallimit}), $P^{\mathrm{CLAIC}}_1\to0$ and
$P^{\mathrm{CLBIC}}_1\to C_3\in(0,1)$.
\end{theorem}

In the model misspecification cases, it is more difficult to compare
analytically the probabilities of selecting model 1 for AIC and CLAIC.
To compare AIC and CLAIC, simulation examples are provided in
Section~\ref{sec-simulation}.\vspace*{-3pt}

%s4 #&#
\section{Simulation studies}
\label{sec-simulation}

In this section, we show simulation results of the following
comparisons in their
decisions among competing models,
\begin{enumerate}
\item CLAIC versus CLBIC,
\item CLAIC versus AIC,
\item CLBIC versus BIC.
\end{enumerate}
To do this, we choose models where the maximum likelihood estimators are
also computationally feasible.
The analysis is different from that in Gao and Song \cite{GS} in that our
concern is not in whether the correct model is asymptotically chosen
with probability 1. If models being compared are close to each other, then
any of the models could be chosen with positive probability, and we are
interested in where CLAIC and AIC might differ.

One general model that allows a variety of univariate and dependence parameters
is the mixed-effects model (see Laird and Ware \cite{LW}); it is
defined via:\vspace*{-1pt}
\begin{eqnarray*}
\Ybf_i&=&\xbf_i\bebf+\mathbf{z}_i
\bbf_i+\epsbf_i,\qquad i=1,2,\ldots,n,
\\
\bbf_i&\sim& N(\zero,\Psi), \qquad\epsbf_i\sim N(\zero,\phi
\Imat _d),
\end{eqnarray*}
where $\bebf$ is $(s+1)$-dimensional vector of fixed effects, $\bbf
_i$ is
$r$-dimensional vector of random effects. $\xbf_i$ and $\mathbf{z}_i$ are
$d\times (s+1)$ and $d\times r$ observable matrices, $\xbf_i$
has a first\vadjust{\goodbreak} column of 1s,
$\phi$ is a variance parameter, $\Psi$ is a $r\times r$
covariance
matrix. Both full likelihood and composite likelihood of the
mixed-effects model can be expressed explicitly with the matrix algebra
notation (see, e.g., Fackler \cite{Fa}, Magnus and Neudecker
\cite{MN}). This model leads to
closed form expressions where $\H$ and $\J$ can be computed (see
Appendix \ref{appb}).

A special case is the clustered data model with exchangeable dependence
structure. It is defined by setting $\mathbf{z}_i=(1,1,\ldots,1)^T$,
$\Psi
=\sigma^2\rho$, and $\phi=\sigma^2(1-\rho)$,
and closed forms for $\H$ and $\J$ can be found in Joe and Lee \cite{JL}.

The three examples given below are representative cases to show patterns
in the decisions from various criteria and in the penalty term
$\tr(\J\H^{-1})$; the patterns were seen over different parameter settings
and dimension $d$.
In the following examples, the composite likelihood corresponding to
the pairwise likelihood or bivariate composite likelihood (BCL) is
specified via\vspace*{-1pt}
\[
S_q=\bigl\{(i,j) \mbox{ for all } i<j\bigr\}.
\]
In Example \ref{ex2}, trivariate composite likelihood (TCL) is also used. The
sets $S_q$ for defining TCL are
\[
S_q=\bigl\{(i,j,k) \mbox{ for all } i<j<k\bigr\}.\vadjust{\goodbreak}
\]
In order that decisions based on AIC and CLAIC are not always for one
model, parameters are
chosen appropriately so that the simpler model has some chance to be chosen.
In Example \ref{ex1}, we consider smaller beta versus larger beta values.

\begin{example}[(Cluster model with exchangeable
covariance matrix,
regression vector $\bebf$ at varying distance from $\zero$)]\label{ex1}
The true number of covariates is 3. Let
$\bebf_0=(0.3,1.3,0.00,0.00)$,
$\bebf_1=(0.3,1.3,0.05,0.02)$,
$\bebf_2=(0.3,1.3,0.15,0.05)$, and
$\bebf_3=(0.3,1.3,0.15,0.10)$,
with first element of the $\bebf$ vectors being the intercept.
Because the last two parameters (regression coefficients for second and
third covariates are smaller), for model selection, simpler models without
the additional covariates might be chosen for any information criteria.
The parameters $\sigma^2=1$ and $\rho=0.5$ are fixed.

%t1 #&#
\begin{table}
\caption{Comparison of decisions for AIC versus CLAIC for different
$\bebf$ vectors,
and distribution of $\tr(\J\H^{-1})$. Cluster size $d=4$;
$\bebf_0=(0.3,1.3,0.00,0.00)$,
$\bebf_1=(0.3,1.3,0.05,0.02)$,
$\bebf_2=(0.3,1.3,0.15,0.05)$,
$\bebf_3=(0.3,1.3,0.15,0.10)$. Covariates are drawn from $N(0,\Imat)$}
\label{tab1}
\begin{tabular*}{\textwidth}{@{\extracolsep{4in minus 4in}}ld{3.0}d{3.0}cd{3.0}d{3.0}d{3.0}cd{3.0}d{3.0}cd{3.0}c@{}}
\hline
& \multicolumn{3}{l}{$\bebf_0$} & \multicolumn{3}{l}{$\bebf_1$} &
\multicolumn{3}{l}{$\bebf_2$} &
\multicolumn{3}{l@{}}{$\bebf_3$} \\[-6pt]
& \multicolumn{3}{l}{\hrulefill} & \multicolumn{3}{l}{\hrulefill} &
\multicolumn{3}{l}{\hrulefill} &
\multicolumn{3}{l@{}}{\hrulefill} \\
\multicolumn{1}{@{}l}{CLAIC\bksl AIC} &
\multicolumn{1}{l}{$1$} &
\multicolumn{1}{l}{$2$} &
\multicolumn{1}{l}{$3$} &
\multicolumn{1}{l}{$1$} &
\multicolumn{1}{l}{$2$} &
\multicolumn{1}{l}{$3$} &
\multicolumn{1}{l}{$1$} &
\multicolumn{1}{l}{$2$} &
\multicolumn{1}{l}{$3$} &
\multicolumn{1}{l}{$1$} &
\multicolumn{1}{l}{$2$} &
\multicolumn{1}{l@{}}{$3$} \\
\hline
$n_c=1$ & 472 & 39 & 11 & 439 & 37 & 19 & 5 & 4 & 1 & 1 & 0 & \phantom{00}1 \\
$n_c=2$ & 23 & 335 & 11 & 24 & 309 & 16 & 1 & 500 & 50 & 1 & 112 & \phantom{0}25 \\
$n_c=3$ & 16 & 5 & 88 & 6 & 13 & 137 & 0 & 33 & 406 & 0 & 13 & 847 \\
\hline
\end{tabular*}
\begin{tabular*}{\textwidth}{@{\extracolsep{\fill}}lcccccccc@{}}
\hline
& \multicolumn{8}{l}{Lower quartile Q1 to upper quartile Q3 of $\tr
(\J\H^{-1})$} \\[-6pt]
& \multicolumn{8}{l@{}}{\hrulefill} \\
& \multicolumn{2}{l}{$\bebf_0$}  & \multicolumn{2}{l}{$\bebf_1$} &
 \multicolumn{2}{l}{$\bebf_2$} &
\multicolumn{2}{l}{$\bebf_3$}  \\[-6pt]
& \multicolumn{2}{l}{\hrulefill} & \multicolumn{2}{l}{\hrulefill} &
\multicolumn{2}{l}{\hrulefill} &
\multicolumn{2}{l@{}}{\hrulefill} \\
\multicolumn{1}{@{}l}{\#covariates} &
\multicolumn{1}{l}{Q1} &
\multicolumn{1}{l}{Q3}  &
\multicolumn{1}{l}{Q1} &
\multicolumn{1}{l}{Q3} &
\multicolumn{1}{l}{Q1} &
\multicolumn{1}{l}{Q3} &
\multicolumn{1}{l}{Q1} &
\multicolumn{1}{l@{}}{Q3}  \\
\hline
1 & 13.7 & 14.1 & 13.7 & 14.1 & 13.6 & 14.0 & 13.6 & 14.0  \\
2 & 16.4 & 16.7 & 16.4 & 16.7 & 16.4 & 16.7 & 16.3 & 16.7  \\
3 & 19.1 & 19.3 & 19.1 & 19.3 & 19.1 & 19.3 & 19.1 & 19.3  \\
\hline
\end{tabular*}
\end{table}

%t2 #&#
\begin{table}[b]
\caption{Comparison of decisions for AIC versus CLAIC for different
$\bebf$ vectors,
and distribution of $\tr(\J\H^{-1})$. Cluster size $d=4$;
$\bebf_0=(0.3,1.3,0.00,0.00)$,
$\bebf_1=(0.3,1.3,0.05,0.02)$,
$\bebf_2=(0.3,1.3,0.15,0.05)$,
$\bebf_3=(0.3,1.3,0.15,0.10)$. Covariates are drawn from $N(0,0.2\Imat
+0.8 \mathbf{1}\mathbf{1}^T)$}
\label{tab2}
\begin{tabular*}{\textwidth}{@{\extracolsep{4in minus 4in}}ld{3.0}d{3.0}cd{3.0}d{3.0}d{3.0}cd{3.0}d{3.0}cd{3.0}c@{}}
\hline
& \multicolumn{3}{l}{$\bebf_0$} & \multicolumn{3}{l}{$\bebf_1$} &
\multicolumn{3}{l}{$\bebf_2$} &
\multicolumn{3}{l@{}}{$\bebf_3$} \\[-6pt]
& \multicolumn{3}{l}{\hrulefill} & \multicolumn{3}{l}{\hrulefill} &
\multicolumn{3}{l}{\hrulefill} &
\multicolumn{3}{l@{}}{\hrulefill} \\
\multicolumn{1}{@{}l}{CLAIC\bksl AIC} &
\multicolumn{1}{l}{$1$} &
\multicolumn{1}{l}{$2$} &
\multicolumn{1}{l}{$3$} &
\multicolumn{1}{l}{$1$} &
\multicolumn{1}{l}{$2$} &
\multicolumn{1}{l}{$3$} &
\multicolumn{1}{l}{$1$} &
\multicolumn{1}{l}{$2$} &
\multicolumn{1}{l}{$3$} &
\multicolumn{1}{l}{$1$} &
\multicolumn{1}{l}{$2$} &
\multicolumn{1}{l@{}}{$3$} \\
\hline
$n_c=1$ & 608 & 43 & 16 & 567 & 41 & 18 & 77 & 24 & 5 & 33 & 9 & \phantom{00}5 \\
$n_c=2$ & 28 & 205 & 7 & 18 & 234 & 8 & 2 & 617 & 37 & 3 & 441 & \phantom{0}39\\
$n_c=3$ & 15 & 2 & 76 & 11 & 3 & 100 & 3 & 28 & 207 & 1 & 29 & 440\\
\hline
\end{tabular*}
\begin{tabular*}{\textwidth}{@{\extracolsep{\fill}}lcccccccc@{}}
\hline
& \multicolumn{8}{l}{Lower quartile Q1 to upper quartile Q3 of $\tr
(\J\H^{-1})$} \\[-6pt]
& \multicolumn{8}{l@{}}{\hrulefill} \\
& \multicolumn{2}{l}{$\bebf_0$}  & \multicolumn{2}{l}{$\bebf_1$} &
 \multicolumn{2}{l}{$\bebf_2$} &
\multicolumn{2}{l}{$\bebf_3$}  \\[-6pt]
& \multicolumn{2}{l}{\hrulefill} & \multicolumn{2}{l}{\hrulefill} &
\multicolumn{2}{l}{\hrulefill} &
\multicolumn{2}{l@{}}{\hrulefill} \\
\multicolumn{1}{@{}l}{\#covariates} &
\multicolumn{1}{l}{Q1} &
\multicolumn{1}{l}{Q3}  &
\multicolumn{1}{l}{Q1} &
\multicolumn{1}{l}{Q3} &
\multicolumn{1}{l}{Q1} &
\multicolumn{1}{l}{Q3} &
\multicolumn{1}{l}{Q1} &
\multicolumn{1}{l@{}}{Q3}  \\
\hline
1 & 13.7 & 14.1 & 13.7 & 14.1 & 13.6 & 14.0 & 13.6 & 14.0  \\
2 & 16.4 & 16.7 & 16.4 & 16.7 & 16.4 & 16.7 & 16.4 & 16.7  \\
3 & 19.1 & 19.3 & 19.1 & 19.3 & 19.1 & 19.3 & 19.1 & 19.3  \\
\hline
\end{tabular*}
\end{table}

%t3 #&#
\begin{table}
\caption{Comparison of decisions for AIC versus CLAIC for different
$\bebf$ vectors,
and distribution of $\tr(\J\H^{-1})$. Cluster size $d=4$;
$\bebf_0=(0.3,1.3,0.00,0.00)$,
$\bebf_1=(0.3,1.3,0.05,0.02)$,
$\bebf_2=(0.3,1.3,0.15,0.05)$,
$\bebf_3=(0.3,1.3,0.15,0.10)$. Covariates are drawn from multivariate
$t$ distribution with $\Sigma_X=\Imat$}
\label{tab5}
\begin{tabular*}{\textwidth}{@{\extracolsep{4in minus 4in}}ld{3.0}d{3.0}cd{3.0}d{3.0}d{3.0}cd{3.0}d{3.0}cd{3.0}c@{}}
\hline
& \multicolumn{3}{l}{$\bebf_0$} & \multicolumn{3}{l}{$\bebf_1$} &
\multicolumn{3}{l}{$\bebf_2$} &
\multicolumn{3}{l@{}}{$\bebf_3$} \\[-6pt]
& \multicolumn{3}{l}{\hrulefill} & \multicolumn{3}{l}{\hrulefill} &
\multicolumn{3}{l}{\hrulefill} &
\multicolumn{3}{l@{}}{\hrulefill} \\
\multicolumn{1}{@{}l}{CLAIC\bksl AIC} &
\multicolumn{1}{l}{$1$} &
\multicolumn{1}{l}{$2$} &
\multicolumn{1}{l}{$3$} &
\multicolumn{1}{l}{$1$} &
\multicolumn{1}{l}{$2$} &
\multicolumn{1}{l}{$3$} &
\multicolumn{1}{l}{$1$} &
\multicolumn{1}{l}{$2$} &
\multicolumn{1}{l}{$3$} &
\multicolumn{1}{l}{$1$} &
\multicolumn{1}{l}{$2$} &
\multicolumn{1}{l@{}}{$3$} \\
\hline
$n_c=1$ & 484 & 43 & 8 & 449 & 42 & 8 & 3 & 5 & 1 & 1 & 0 & \phantom{00}1 \\
$n_c=2$ & 28 & 314 & 13 & 21 & 311 & 17 & 1 & 539 & 43 & 1 & 129 & \phantom{0}33 \\
$n_c=3$ & 12 & 16 & 82 & 18 & 9 & 125 & 0 & 32 & 376 & 0 & 17 & 818 \\
\hline
\end{tabular*}
\begin{tabular*}{\textwidth}{@{\extracolsep{\fill}}lcccccccc@{}}
\hline
& \multicolumn{8}{l}{Lower quartile Q1 to upper quartile Q3 of $\tr
(\J\H^{-1})$} \\[-6pt]
& \multicolumn{8}{l@{}}{\hrulefill} \\
& \multicolumn{2}{l}{$\bebf_0$}  & \multicolumn{2}{l}{$\bebf_1$} &
 \multicolumn{2}{l}{$\bebf_2$} &
\multicolumn{2}{l}{$\bebf_3$}  \\[-6pt]
& \multicolumn{2}{l}{\hrulefill} & \multicolumn{2}{l}{\hrulefill} &
\multicolumn{2}{l}{\hrulefill} &
\multicolumn{2}{l@{}}{\hrulefill} \\
\multicolumn{1}{@{}l}{\#covariates} &
\multicolumn{1}{l}{Q1} &
\multicolumn{1}{l}{Q3}  &
\multicolumn{1}{l}{Q1} &
\multicolumn{1}{l}{Q3} &
\multicolumn{1}{l}{Q1} &
\multicolumn{1}{l}{Q3} &
\multicolumn{1}{l}{Q1} &
\multicolumn{1}{l@{}}{Q3}  \\
\hline
1 & 14.8 & 15.2 & 14.8 & 15.2 & 14.7 & 15.1 & 14.7 & 15.1 \\
2 & 17.2 & 17.5 & 17.2 & 17.5 & 17.2 & 17.5 & 17.2 & 17.5 \\
3 & 19.7 & 19.8 & 19.7 & 19.8 & 19.7 & 19.8 & 19.7 & 19.8 \\
\hline
\end{tabular*}
\end{table}

For each of the four $\bebf$ vectors, 1000 replicates with sample size
$n=100$ and cluster size $d=4$ are generated. Three different settings
are used to simulate the covariates and the random effects. In settings
(i) and (ii), the covariates $\xbf_i=(\xbf_{i1},\xbf_{i2},\xbf_{i3})^T$
are independent random vectors from $N(0,\Sigma_X)$ with $\Sigma
_X=\Imat
$, the identity matrix and $\Sigma_X=0.2\Imat+0.8\mathbf{1}\mathbf{1}^T$,
respectively. The random effect $\bbf_i$ is obtained from normal
distribution. In setting (iii), $t$-distribution with degree of freedom
$3$ is used for $\bbf_i$ instead so that the robustness of the
information criteria under model misspecification can therefore be
investigated. That is, $\bbf_i=\rho^{1/2}t_i$, where $t_i$ are
independent $t$-distributed random variables. We then compare the
decisions of AIC and CLAIC for regression models with the first, the
first two or all three covariates ($n_c=1,2$ or 3). For setting (i),
summaries in Table~\ref{tab1} show patterns in the decisions and
in the amount of variation in the CLAIC penalty term $\tr(\J\H^{-1})$.
As an example, for $\bebf_1$, there were 137 cases where both AIC and CLAIC
chose the 3-covariate model.
Table~\ref{tab1} shows that the decisions for CLAIC are the same as
with AIC in a
high proportion of cases;
both tend to choose a regression model with more covariates if
the true $\bebf$ vector has more coefficients farther from 0.
The results of BIC and CLBIC are similar, and are not shown.\vadjust{\goodbreak}
The variation in $\tr(\J\H^{-1})$ is not too much when the sample size
is large enough. As implied by Lemma \ref{mono}, $\tr(\J\H^{-1})$
tends to increase for models
with additional parameters. Similar results of settings (ii) and (iii)
are given in Tables~\ref{tab2} and~\ref{tab5}, respectively. In
this example, it can be seen that all information criteria have a
higher chance to select the smaller model in the presence of strong
correlations (say, 0.8) in the covariates. In the case where the
distribution of $\epsbf_i$ is misspecified, the decisions from all
information criteria are very similar to the counterpart without
misspecification.
\end{example}

\begin{example}[(Multivariate normal regression model, different covariance
structures)]\label{ex2}
This example shows local alternatives or perturbations of different types,
either in univariate or in dependence parameters.
We compare exchangeable (exch) versus unstructured (unstr) dependence
when true
covariance matrix has different deviations from exchangeable.
The choices of the true covariance matrices are:
\begin{eqnarray*}
\Sigma_1&=& %
\pmatrix{1 &0.5 &0.5 &0.5
\cr
0.5 & 1 &0.5 &0.5
\cr
0.5 &0.5 & 1 &0.5
\cr
0.5 &0.5 &0.5 & 1
} %
, \\ %\mbox{and}
\Sigma_2&=& %
\pmatrix{1 &0.5+\varepsilon_1/\sqrt{n}
&0.5 &0.5
\cr
0.5+\varepsilon _1/\sqrt{n} & 1 &0.5 &0.5
\cr
0.5 &0.5 & 1 &
0.5+\varepsilon_1/\sqrt{n}
\cr
0.5 &0.5 &0.5+\varepsilon _1/
\sqrt {n} & 1 },
\\
\Sigma_{ka}&=&\diag(1,1,1+\varepsilon_2/\rt{n},1+
\varepsilon_2/\rt{n}) \Sigma _k \diag(1,1,1+
\varepsilon_2/\rt{n},1+\varepsilon_2/\rt{n})
\end{eqnarray*}
for $k=1,2$,
where $\varepsilon_1=0.07\sqrt{200}$ and $\varepsilon_2=0.05\sqrt{200}$.
$\Sigma_2$ changes some correlation parameters,
$\Sigma_{1a}$ changes some variance (univariate) parameters, and
$\Sigma_{2a}$ changes both correlation and variance parameters.
The regression vector $\bebf=(0.3,1.3)$ is fixed and the covariates
$\xbf_i$ are independent standard normal random variables. Summaries in
Table~\ref{tab3} are from 1000 replicates with different sample
sizes $n$ and cluster size $d=4$.

%t4 #&#
\begin{table}
\caption{Comparison of decisions for AIC versus CLAIC under different
perturbations
of the exchangeable dependence model}
\label{tab3}
\begin{tabular*}{\textwidth}{@{\extracolsep{\fill}}lllllllll@{}}
\hline
& \multicolumn{2}{l}{$\Sigma_1$} & \multicolumn{2}{l}{$\Sigma
_{1a}$} &
\multicolumn{2}{l}{$\Sigma_{2}$} & \multicolumn{2}{l@{}}{$\Sigma
_{2a}$} \\[-6pt]
& \multicolumn{2}{l}{\hrulefill} & \multicolumn{2}{l}{\hrulefill} &
\multicolumn{2}{l}{\hrulefill} & \multicolumn{2}{l@{}}{\hrulefill} \\
CLAIC\bksl AIC & exch. & unstr. & exch. & unstr. & exch. & unstr. & exch. &
unstr. \\
\hline
\multicolumn{9}{@{}l}{$n=200$, $d=4$, BCL} \\
exch. & 919 & \phantom{0}16 & 813 & \phantom{0}15 & 668 & 211 & 574 & 162 \\
unstr. & \phantom{0}40 & \phantom{0}25 & \phantom{0}95 & \phantom{0}77 & \phantom{0}21 & 100 & \phantom{0}45 & 219 \\[3pt]
\multicolumn{9}{@{}l}{$n=500$, $d=4$, BCL} \\
%CLAIC\bksl AIC & exch & unstr & exch & unstr & exch & unstr & exch &
%unstr \\
exch. & 911 & \phantom{0}12 & 825 & \phantom{0}10 & 699 & 175 & 593 & 168 \\
unstr. & \phantom{0}50 & \phantom{0}27 & 108 & \phantom{0}57 & \phantom{0}18 & 108 & \phantom{0}42 & 197 \\[3pt]
\multicolumn{9}{@{}l}{$n=500$, $d=4$, TCL} \\
%CLAIC\bksl AIC & exch & unstr & exch & unstr & exch & unstr & exch &
%unstr \\
exch. & 944 & \phantom{00}6 & 890 & \phantom{00}6 & 710 & \phantom{0}94 & 617 & \phantom{0}84 \\
unstr. & \phantom{0}17 & \phantom{0}33 & \phantom{0}43 & \phantom{0}61 & \phantom{00}7 & 189 & \phantom{0}18 & 281 \\
\hline
\end{tabular*}
\end{table}

The patterns are similar to above for larger cluster size
$d=5,6,7$ and perturbations of a different exchangeable correlation matrix.
That is, CLAIC tends to more often than AIC choose the unstructured
dependence when the perturbation is only in the variances (i.e.,
$\Sigma_{1a}$), and AIC tends to more often than CLAIC choose the
unstructured dependence when the perturbation is only in the correlations
(i.e., $\Sigma_2$). For perturbations in the correlations, going to
trivariate composite likelihood
makes CLAIC closer to AIC in the decision between the two models.

For $\Sigma_1$,
CLAIC selects bigger model more often than AIC in all three settings
(see Table~\ref{tab3}). However, the probabilities $P^{\mathrm{CLAIC}}_1$ and
$P^{\mathrm{AIC}}_1$ are very close to each other. In this example, CLAIC and
AIC give very similar decisions. The outcome is consistent with Theorem
\ref{main}(2). Under $H_1$, AIC selects model 1 with probability
approximately $\Pr(Z^2_1+\cdots+Z^2_m<2m)$. For the TCL with $S_q=\{
(i,j,k)\}$, $n=500$, $\Sigma=\Sigma_1$, CLAIC selects model 1 with
probability approximately $\Pr(\lambda_1Z^2_1+\cdots+\lambda
_mZ^2_m<2(\lambda_1+\cdots+\lambda_m))$. Here, $\lambda_1,\ldots
,\lambda_m$ are
\[
3.34, 2.87, 2.73, 2.52, 2.07, 2.03, 1.61, 1.50.
\]
Since these $\lambda$ values differ from each other, Lemma \ref{prob}(2) guarantees that
\[
\Pr\bigl(\lambda_1Z^2_1+\cdots+
\lambda_mZ^2_m<2(\lambda_1+\cdots
+\lambda _m)\bigr)<\Pr\bigl(Z^2_1+
\cdots+Z^2_m<2m\bigr).
\]
Indeed, for the eigenvalues $\lambda$ in the example, we have
\[
\Pr\bigl(\lambda_1Z^2_1+\cdots+
\lambda_8Z^2_8>2(\lambda_1+
\cdots +\lambda _8)\bigr)=0.0468 \quad\mbox{and}\quad \Pr\bigl(Z^2_1+
\cdots+Z^2_8>16\bigr)=0.0424.
\]
Here, the numerical method proposed in Rice \cite{Ri} is used to
obtain the
first probability. The first probability is slightly greater than the
second probability.
\end{example}

\begin{example}[(Multivariate normal regression model, different covariance structures)]\label{ex3}
This example shows the exchangeable (exch) dependence model and its local
alternatives with perturbations of different sizes in dependence
parameters. Information criteria AIC, BIC, CLAIC, and CLBIC are
compared. The choices of the true covariance matrices are:
\[
\Sigma(\delta)= %
\pmatrix{1 &0.5+\delta&0.5 &0.5
\cr
0.5+\delta& 1 &0.5
&0.5
\cr
0.5 &0.5 & 1 &0.5+\delta
\cr
0.5 &0.5 &0.5+\delta& 1} .
\]
Define $\Sigma_1=\Sigma(0)$, $\Sigma_2=\Sigma(n^{-1/2})$,
$\Sigma
_3=\Sigma(n^{-1/2}\log n)$, and $\Sigma_4=\Sigma(0.2)$.
The regression vector $\bebf=(0.3,1.3)$ is fixed and the covariates
$\xbf_i$ are independent standard normal random variables.

%t5 #&#
\begin{table}
\tablewidth=300pt
\caption{Comparison of decisions for AIC, BIC, CLAIC, and CLBIC under
different perturbations
of the exchangeable dependence model. Sample size $n=500$}
\label{tab4}
\begin{tabular*}{300pt}{@{\extracolsep{4in minus 4in}}lllll@{}}
\hline
&\multicolumn{4}{l}{Frequency of selecting exchangeable}\\[-6pt]
&\multicolumn{4}{l@{}}{\hrulefill}\\
Info. crit. & $\Sigma_1$ & $\Sigma_2$ & $\Sigma_3$ & \multicolumn{1}{l@{}}{$\Sigma_4$}\\
\hline
AIC & \phantom{0}961 & \phantom{0}712 & \phantom{00}5 & 0 \\
CLAIC & \phantom{0}950 & \phantom{0}803 & \phantom{0}28 & 0 \\
BIC & 1000 & 1000 & 705 & 0 \\
CLBIC & 1000 & 1000 & 927 & 0 \\
\hline
\end{tabular*}
\end{table}

Summaries in Table~\ref{tab4} are from 1000 replicates sample
size $n=500$ and cluster size  \mbox{$d=4$}. The frequencies of selecting the
exchangeable dependence model are reported.
We see that BIC/CLBIC tends to select the exchangeable dependence model
more often than AIC/CLAIC. Under the assumption of exchangeable
dependence model, BIC/CLBIC have greater chance of selecting the
correct model. However, BIC/CLBIC are less sensitive to small
perturbations than AIC/CLAIC. The results are consistent with Theorem
\ref{main}.\vadjust{\goodbreak}
\end{example}

\begin{example}[(Comparison of information criteria under
model misspecification)]\label{ex4}
To see the effect of model misspecification, we
repeat Example \ref{ex3} with the following changes: (i) $\bbf_i=\mathbf
{Cu}_i$,
where $\Psi=\mathbf{CC}^T$ is the Cholesky decomposition and $\mathbf{u}_i$
are vectors of independent Laplace random variables with mean zero and
variance one.
(ii) $\bbf_i$ is generated from normal distribution but
\[
\Ybf_i=(0.3,0.6,0.9,1.2)^T+\xbf_i\bebf+
\mathbf{z}_i \bbf_i+\epsbf _i,\qquad i=1,2,
\ldots,n.
\]
The results under (i) and (ii) are summarized in Tables~\ref{tab6} and
\ref{tab7}, respectively. The decisions under (i) are comparable to
(1), (3), (4) in Theorem \ref{main.robust}. The decisions under (ii)
are similar to that described in (2) in Theorem \ref{main.robust}.
Comparing with Example \ref{ex3}, under both (i) and (ii), the alternative
model is more likely to be selected.
\end{example}

%t6 #&#
\begin{table}
\tablewidth=300pt
\caption{Comparison of decisions for AIC, BIC, CLAIC, and CLBIC under
perturbation
in the distribution law. Sample size $n=500$}
\label{tab6}
\begin{tabular*}{300pt}{@{\extracolsep{4in minus 4in}}lllll@{}}
\hline
&\multicolumn{4}{l}{Frequency of selecting exchangeable}\\[-6pt]
&\multicolumn{4}{l@{}}{\hrulefill}\\
Info. crit. & $\Sigma_1$ & $\Sigma_2$ & $\Sigma_3$ & \multicolumn{1}{l@{}}{$\Sigma_4$}\\
\hline
AIC & \phantom{0}768 & \phantom{0}479 & \phantom{00}3 & 0 \\
CLAIC & \phantom{0}712 & \phantom{0}537 & \phantom{0}12 & 0 \\
BIC & 1000 & \phantom{0}998 & 575 & 0 \\
CLBIC & 1000 & 1000 & 795 & 0 \\
\hline
\end{tabular*}
\end{table}

%t7 #&#
\begin{table}
\tablewidth=300pt
\caption{Comparison of decisions for AIC, BIC, CLAIC, and CLBIC under
perturbation
in the mean. Sample size $n=500$}
\label{tab7}
\begin{tabular*}{300pt}{@{\extracolsep{4in minus 4in}}lllll@{}}
\hline
&\multicolumn{4}{l}{Frequency of selecting exchangeable}\\[-6pt]
&\multicolumn{4}{l@{}}{\hrulefill}\\
Info. crit. & $\Sigma_1$ & $\Sigma_2$ & $\Sigma_3$ & \multicolumn{1}{l@{}}{$\Sigma_4$}\\
\hline
AIC & \phantom{0}0 & 0 & 0 & 0 \\
CLAIC & \phantom{0}0 & 0 & 0 & 0 \\
BIC & 24 & 0 & 0 & 0 \\
CLBIC & 68 & 3 & 0 & 0 \\
\hline
\end{tabular*}
\end{table}

%s5 #&#
\section{Spruce tree growth data}
\label{sec-realdata}

In this section, we study the spruce tree growth data in Example 1.3 in
Diggle \textit{et al.} \cite{DLZ}. The decisions from AIC (BIC) and their
composite likelihood
counterparts are compared.

The dataset consists of the data from $n=79$ trees and is available
in the R package MEMSS (Pinheiro and Bates \cite{PB}).
For each tree, the logarithm of the volume of the tree trunk was
estimated and recorded in $d=13$ chosen days $t_1,t_2,\ldots,t_{13}$
from the beginning of the experiment over a period of 674 days. The
trees were grown in four different plots, labeled $1,2,3,4$, respectively.
The days are 152, 174, 201, 227, 258, 469, 496, 528, 556, 579, 613,
639, 674
days since 1988-01-01, corresponding to roughly beginning of June to
mid-August in 1988 and mid-April to the end of October in 1989.
The first two plots represent an ozone-controlled atmosphere
and the last two plots represent a normal atmosphere.
From the plots in Diggle \textit{et al.} \cite{DLZ}, the growth rates in the
two time periods
are different.

A linear mixed-effects model accounts for different growth
rates in the two periods is the following.
For a given tree, with $y=\log$ size has growth and $t=\mathrm{day}$ since 1988-01-01,
\begin{eqnarray*}
y_{ij}&=&a_0+a_1(t_j-152)/100+
\varepsilon_i, \qquad 152\le t_j \le258,
\\
y_{ij}&=&\bigl[a_0+a_1(258-152)/100\bigr]+
a_2(t_j-445)/100+\varepsilon_i,\qquad 469\le
t_j \le674.
\end{eqnarray*}
To introduce fixed and random effects, $a_0=\b_0+\b_3I(\mathrm{ozone})+b_0$,
where $b_0$ is random with normal distribution; in addition,
$a_1=\b_1+\b_4I(\mathrm{ozone})+b_1$, $a_2=\b_2+\b_5I(\mathrm{ozone})+b_2$,
where $b_1,b_2$ are also random and normally distributed.
There was little growth in between the two periods so the use of $445=469-24$
treats the days 258 and 469 as one measurement unit apart.

%t8 #&#
\begin{table}
\caption{Spruce data: Comparison of parameter estimates from maximizing
full likelihood, TCL, BCL; the correlation of the random effects are
small and not included.
Standard errors (SEs) are obtained via the delete-one jackkknife}
\label{tab8}
\begin{tabular*}{\textwidth}{@{\extracolsep{4in minus 4in}}ld{2.10}d{2.10}c@{}}
\hline
\multicolumn{1}{@{}}{Parameter} & \multicolumn{1}{l}{Full (SE)} & \multicolumn{1}{l}{TCL (SE)} & \multicolumn{1}{l@{}}{BCL (SE)} \\
\hline
$\b_0$ & 4.272\ (0.154) & 4.310\ (0.152) & \phantom{$-$}4.311 (0.152) \\
$\b_1$ & 1.415\ (0.064) & 1.371\ (0.062) & \phantom{$-$}1.373 (0.062) \\
$\b_2$ & 0.371\ (0.021) & 0.383\ (0.021) & \phantom{$-$}0.382 (0.021) \\
$\b_3$ & -0.101\ (0.173) & -0.097\ (0.171) & $-$0.097 (0.171) \\
$\b_4$ & -0.223\ (0.076) & -0.228\ (0.074) & $-$0.227 (0.075) \\
$\b_5$ & -0.012\ (0.027) & -0.012\ (0.027) & $-$0.012 (0.027) \\[3pt]
$\operatorname{resid}$SD & 0.138\ (0.005) & 0.126 (0.005) & \phantom{$-$}0.118 (0.006) \\
SD($b_0$) & 0.616\ (0.051) & 0.625\ (0.050) & \phantom{$-$}0.630 (0.050) \\
SD($b_1$) & 0.270\ (0.031) & 0.323\ (0.030) & \phantom{$-$}0.353 (0.034) \\
SD($b_2$) & 0.098\ (0.018) & 0.110\ (0.017) & \phantom{$-$}0.118 (0.017) \\
\hline
\end{tabular*}
\end{table}

Estimates of regression coefficients for the fixed effects and SDs
of the random effects are shown in Table~\ref{tab8}; the standard errors
of these parameter estimates are obtained with the delete-one jackknife
as mentioned in Varin \textit{et al.} \cite{VRF} for composite likelihood methods.
Based on the estimates in this table, for submodels we consider setting
$\b_5,\b_3,\b_4$ in turn to zero for the effects of ozone in the
second period, initial point, and first period.
Hence, we have submodels with 5, 4 and 3 regression parameters.
In Table~\ref{tab9}, the decisions of the difference full
likelihood and composite likelihood information criteria are shown.

For these four models, all of the information criteria chose the same
best model with a significant $\b_4$, the effect of ozone for the growth
rate in the first period. Based on these criteria and standard errors,
the effect $\b_5$ of the ozone for the growth in the
second period is much more negligible, and the effect $\b_3$ of ozone
for the period before day 152 is also non-significant.
Note that the model with $\b_5=0$ and five non-zero $\b$'s,
the AIC/BIC values are relatively closer
to those for the best model than the corresponding CLAIC/CLBIC values;
this is also seen in the corresponding $z$-statistics:
for $\b_3$, the ratio of estimate and SE
is $-0.118/0.162=-0.73$ for full likelihood,
$-0.097/0.171=-0.57$ for TCL,
and $-0.094/0.175=-0.54$ for BCL.

Although the four models in Table~\ref{tab9} are ranked the same
on all information criteria, this is not the case when we also
consider other models with additional binary variables
to handle four plots (two plots for each of ozone and control).
That is, to relate to what we found in the simulation examples
in Section~\ref{sec-simulation}, if we consider many models and some
of them are quite close in fit because of some regression coefficients
being near zero, then the rankings can be different for
full and composite likelihood information criteria.

%t9 #&#
\begin{table}
\caption{Spruce data: Comparison of decisions for AIC, BIC, CLAIC, and CLBIC.
The decision is the number of
$\b$'s in the model with smallest information criterion value.
The values of CLAIC and CLBIC have been divided by
${13\choose3}=286$ for TCL and by ${13\choose2}=78$ for BCL in order
that they are smaller}
\label{tab9}
\begin{tabular*}{\textwidth}{@{\extracolsep{\fill}}lllllll@{}}
\hline
& \multicolumn{2}{l}{Full likelihood} & \multicolumn{2}{l}{TCL} &
\multicolumn{2}{l@{}}{BCL}\\
& \multicolumn{2}{l}{\hrulefill} & \multicolumn{2}{l}{\hrulefill} &
\multicolumn{2}{l@{}}{\hrulefill}\\
\#$\b$'s & AIC & BIC & CLAIC & CLBIC & CLAIC & CLBIC\\
\hline
6 & $-2319.5$ & $-2288.7$ & $-291.9$ & $-276.7$ & $-124.0$ & $-112.$1 \\
5 ($\b_5=0$) & $-2321.3$ & $-2292.9$ & $-292.6$ & $-278.3$ & $-124.4$ & $-113.1$ \\
4 ($\b_5=\b_3=0$) & $-2322.7$ & $-2296.6$ & $-293.8$ & $-281.4$ & $-125.4$ &
$-115.6$ \\
3 ($\b_5=\b_3=\b_4=0$) & $-2314.7$ & $-2291.0$ & $-288.4$ & $-277.5$ &
$-120.9$ &
$-112.7$ \\[3pt]
Decision & 4 & 4 & 4 & 4 & 4 & 4\\
\hline
\end{tabular*}
\end{table}

%s6 #&#
\section{Discussion}
\label{sec-discussion}

In this paper, we have results that show how decisions from CLAIC
compare with
those from AIC for nested models.
This was mostly based on the theory of local alternatives applied
to composite likelihood; this is the theory that is most relevant to
understand how model selection performs for models that are not
far apart.

The theory of this paper can be applied to other models
to understand better how CLAIC compares with AIC for different types
of perturbations that may involve univariate or dependence parameters.
This can be done if the $\J$ and $\H$ can be computed, possibly
based on simulation methods.
Further analysis will help in the understanding of conditions for
which CLAIC has more ``power'' to detect a more complex model.
The results have some analogies with those
in Joe and Maydeu-Olivares \cite{JM}, where it is shown that there are
directions of local alternatives for which goodness-of-fit statistics based
on low-dimensional margins can have more power.

Although analysis in this paper is with composite marginal likelihood,
we expect many of the results apply to composite conditional likelihood.

Another topic of research is further study of the extension of the
procedure of Vuong \cite{Vu} for composite likelihood to understand its
potential usefulness for comparing prediction similarity for non-nested
models.

\begin{appendix}\label{app}

%s7 #&#
\section{Proofs}\label{appa1}

%s7.1 #&#
\subsection{Assumptions}\label{appa}

The following assumptions are used, similar to Vuong \cite{Vu}.

 A1: $\Theta$, $\Gamma$ are compact subsets of a Euclidean space.

 A2: Let $\vthbf=\thbf$ for model 1 and $\vthbf=\gbf$ for
model 2.
For $M=1,2$,
under the true model, we have almost surely for all $(\xbf,\mathbf{y})$,
$\log f^{(M)}_{S_q}(\mathbf{y}_{S_q};\xbf,\thbf)$ is twice continuously\vspace*{-1.5pt}
differentiable over the parameter space. In addition, there exist integrable
(under the true model) functions $K^{(M)}_q(\xbf,\mathbf{y})$,
$K^{(M)}_{qj}(\xbf,\mathbf{y})$, $K^{(M)}_{qjk}(\xbf,\mathbf{y})$,
where $\vth_j,\vth_k$ are components in the parameter $\vthbf$, such that
\begin{eqnarray*}
\sup\bigl\llvert \log f^{(M)}_{S_q}(\mathbf{y}_{S_q};
\xbf_i,\vthbf)\bigr\rrvert ^2 &<& K_q^{(M)}(
\xbf,\mathbf{y}),
\\
\sup\biggl\llvert \frac{\p}{\p\vth_j} \log f^{(M)}_{S_q}(\mathbf{y}_{S_q};
\xbf_i,\vthbf)\biggr\rrvert ^2 &<& K^{(M)}_{qj}(
\xbf,\mathbf{y}),
\\
\sup\biggl\llvert \frac{\p^2}{\p\vth_j\,\p\vth_k} \log f^{(M)}_{S_q}(\mathbf{y}_{S_q};
\xbf_i,\vthbf)\biggr\rrvert &<& K^{(M)}_{qjk}(\xbf
,\mathbf{y}),
\end{eqnarray*}
where the suprema are over the parameter space $\Theta$ or $\Gamma$.

 A3: Under the true model, for $f^{(1)}$, the local maximum point
\[
\thbf^*=\arg\max_{\Theta}\lim_{n\to\oo}
n^{-1}\sum_{i=1}^n \E \Biggl\{
\sum_{q=1}^Q\log f^{(1)}_{S_q}(
\Ybf_{i,S_q};\xbf_i,\thbf) \Biggr\}
\]
is unique and $\thbf^*$ is an interior point of $\Theta$.
Similarly $\gbf^*$ is defined for $f^{(2)}$ and is an interior point of
$\Gamma$.

Assumption A2 guarantees the existence of positive definite matrices
$\H^{(1)},\H^{(2)},\J$ given below.
For the matrices defined below, all expectations below are taken under
the true model.
\begin{eqnarray*}
\H^{(1)}(\thbf) &=& -\lim_{n\to\oo}n^{-1}\sum
_{i=1}^n \E_g \Biggl\{\sum
_{q=1}^{Q} \frac{\p^2}{\p\thbf\,\p\thbf^T}\log
f^{(1)}_{{S_q}}(\Ybf _{i,S_q};\xbf _i,\thbf)
\Biggr\},
\\
{\J}^{(11)}(\thbf) &=&\lim_{n\to\oo}n^{-1}\sum
_{i=1}^n \E_g \Biggl\{
\frac{\p}{\p\thbf}\sum_{q=1}^{Q}\log
f^{(1)}_{{S_q}}(\Ybf _{i,S_q};\xbf_i,\thbf)
\cdot\frac{\p}{\p\thbf^T}\sum_{q=1}^{Q}\log
f^{(1)}_{{S_q}}(\Ybf _{i,S_q};\xbf_i,\thbf)
\Biggr\},
\\
{\J}^{(12)}(\thbf,\gbf) &=&\lim_{n\to\oo}n^{-1}
\sum_{i=1}^n \E_g \Biggl\{
\frac{\p}{\p\thbf}\sum_{q=1}^{Q}\log
f^{(1)}_{{S_q}}(\Ybf _{i,S_q};\xbf_i,\thbf)
\cdot\frac{\p}{\p\gbf^T}\sum_{q=1}^{Q}\log
f^{(2)}_{{S_q}}(\Ybf _{i,S_q};\xbf_i,\gbf)
\Biggr\}.
\end{eqnarray*}
Similarly $\H^{(2)}(\gbf)$, ${\J}^{(22)}(\gbf)$,
${\J}^{(21)}(\gbf,\thbf)$ can be defined.
Let
\[
\J=\pmatrix{\J^{(11)}\bigl(\thbf^*
\bigr)&\J^{(12)}\bigl(\thbf^*,\gbf^*\bigr)
\vspace*{2pt}\cr
\J^{(21)}\bigl(\gbf^*,\thbf^*\bigr)&\J^{(22)}\bigl(\gbf^*\bigr)}.
\]
Applying the law of large numbers and the Central Limit theorem, we have
as $n\to\oo$,
%
%e7.1 #&#
%e7.2 #&#
%
\begin{eqnarray}
-n^{-1}\left[ %
\matrix{
\displaystyle\frac{\p^2}{\p\thbf\,\p\thbf^T}L^{(1)}_{\mathrm{CL}}\bigl(\thbf^*\bigr) &
\displaystyle\frac{\p^2}{\p\gbf\,\p\gbf^T}L^{(2)}_{\mathrm{CL}}\bigl(\gbf^*\bigr) }
\right] &\mathop{\longrightarrow}^{a.s.}& \bigl[
\matrix{ \H^{(1)} & \H^{(2)} }\bigr], \label{eq-limhess}
\\
\left[ \matrix{\displaystyle\frac{1}{\sqrt{n}}
\frac{\p}{\p\thbf
^T}L^{(1)}_{\mathrm{CL}}\bigl(\theta^*\bigr)&
\displaystyle\frac{1}{\sqrt{n}}\frac{\p}{\p\gbf^T}L^{(2)}_{\mathrm{CL}}\bigl(\gbf^*
\bigr) }
\right] &\mathop{\longrightarrow}^{d}& N(
\zero,\J). \label{eq-limscore}
\end{eqnarray}

%s7.2 #&#
\subsection{Proof of Proposition \texorpdfstring{\protect\ref{lik-ratio}}{3.1}}\label{appa2}

The proof can be established following the same arguments as in Vuong
\cite{Vu}, so that most details are omitted. Below, the asymptotic
covariance matrix is obtained in a heuristic way.

Based on (\ref{eq-limhess}) and (\ref{eq-limscore}), and the assumptions
A1--A3 (see Appendix \ref{appa1}), Taylor expansions to second order are valid
and lead to:
\[
2\operatorname{LR}= n\bigl(\hat{\gbf}-\gbf^*\bigr)^T\H^{(2)}\bigl(\hat{
\gbf}-\gbf^*\bigr) - n\bigl(\hat{\thbf}-\thbf^*\bigr)^T
\H^{(1)}\bigl(\hat{\thbf}-\thbf^*\bigr) +\mathrm{o}_p(1),
\]
and the matrix of the (asymptotic) quadratic form in independent
standard normal random variables is
$\V^{1/2}\diag(-\H^{(1)},\H^{(2)})\V^{1/2}$,
where
\[
\V=\diag\bigl(\bigl(\H^{(1)}\bigr)^{-1}, \bigl(
\H^{(2)}\bigr)^{-1}\bigr)\J \diag\bigl(\bigl(\H^{(1)}
\bigr)^{-1}, \bigl(\H^{(2)}\bigr)^{-1}\bigr) =\pmatrix{ \V^{(11)}& \V^{(12)}
\vspace*{2pt}\cr
\V^{(21)}& \V^{(22)} }
\]
is the asymptotic covariance matrix of
$n^{1/2}(\thatbf-\thbf^*,\ghatbf-\gbf^*)$.
The eigenvalues of this matrix are the same as those of
\[
\Amat=\pmatrix{ -\H^{(1)}
\V^{(11)}& -\H^{(1)}\V^{(12)}
\vspace*{2pt}\cr
\H^{(2)}\V^{(21)}& \H^{(2)}\V^{(22)}}.
\]
Let
\[
\Kmat=\pmatrix{ -\Imat_{p_1} & \zero
\vspace*{2pt}\cr
\zero& \Imat_{p_2}}.
\]
Then
\[
\Kmat\Amat\Kmat= \pmatrix{ -\bigl(
\J^{(11)}\bigr) \bigl(\H^{(1)}\bigr)^{-1}&\bigl(
\J^{(12)}\bigr) \bigl(\H^{(2)}\bigr)^{-1}
\vspace*{2pt}\cr
-\bigl(\J^{(21)}\bigr) \bigl(\H^{(1)}\bigr)^{-1}&
\bigl(\J^{(22)}\bigr) \bigl(\H^{(2)}\bigr)^{-1}},
\]
and the eigenvalues of this matrix and $\Amat$ are the same.

%s7.3 #&#
\subsection{Proof of Theorem \texorpdfstring{\protect\ref{main}}{3.1}}\label{appa3}

Consider the nested case where
$f^{(1)}(\cdot;\mathbf{x},\thbf)=f^{(2)}(\cdot;\mathbf{x},\thbf,\zero)$.
Suppose that $\gbf=(\thbf,\zebf)$ is $p_2$-dimensional and $\zebf$ is
$m$-dimensional, where
$m=p_2-p_1$.
(Note: the maximum composite likelihood estimator for model 1 is $\hat
{\thbf}$,
and it is not the sub-vector of $\hat{\boldsymbol{\gamma}}$, the maximum composite
likelihood estimator for model 2.)
For convenience, the following notation is used throughout the proof,
\begin{eqnarray*}
\H^{(1)}&=&\H_{\theta\theta} \quad\mbox{and} \quad\H^{(2)}=\pmatrix{ \H_{\theta\theta}&\H_{\theta\zeta}
\vspace*{2pt}\cr
\H_{\zeta\theta}&\H_{\zeta\zeta}},\\
\J^{(11)}&=&\J_{\theta\theta},\qquad \J^{(22)}=\pmatrix{ \J_{\theta\theta}&\J_{\theta\zeta}
\vspace*{2pt}\cr
\J_{\zeta\theta}&\J_{\zeta\zeta} },\qquad
\J^{(21)}=\pmatrix{ \J_{\theta\theta}
\vspace*{2pt}\cr
\J_{\zeta\theta} },\qquad \J^{(12)}=\pmatrix{ \J_{\theta\theta}&\J_{\theta\zeta}
}.
\end{eqnarray*}

\textit{Proof of} (1). For CLBIC, it is a special case of
Theorem 1 and 2 in Gao and Song \cite{GS}. A detailed treatment on the order
consistency can be found in Gao and Song \cite{GS}. Below, we complete
the proof by
showing that $P^{\mathrm{CLAIC}}_1$, the probability that CLAIC selects model
1 under $H_1$ has the form $P(\lambda_1 U_1+\cdots+\lambda_m
U_m<2(\lambda_1+\cdots+\lambda_m))$.

 Let $\gbf^*=(\thbf^*,\zebf^*)$ be the true value. Under the
null hypothesis, $\zebf^*=\zero$.
From Taylor expansions of $L^{(2)}_{\mathrm{CL}}(\gbf^*)$ and
$L^{(1)}_{\mathrm{CL}}(\thbf^*)$
around $\hat{\gbf}$ and $\hat{\thbf}$,
we have the composite log-likelihood ratio:
\[
0\le \operatorname{LR}=\tfrac{1}{2}n\bigl(\hat{\gbf}-\gbf^*\bigr)^T
\H^{(2)}\bigl(\hat{\gbf }-\gbf^*\bigr) -\tfrac{1}{2}n\bigl(\hat{
\thbf}-\thbf^*\bigr)^T\H^{(1)}\bigl(\hat{\thbf}-\thbf^*\bigr)
+\mathrm{o}_p(1).
\]
From Proposition \ref{lik-ratio},
it has asymptotically the same distribution as $\Zbf^TD\Zbf$ where
$\Zbf
$ is a $(p_1+p_2)$-vector of independent $N(0,1)$ random variables and
$D$ is a diagonal matrix with diagonal elements equal to the
eigenvalues of $\Bmat$ (defined in (\ref{eq-Bmat})). In addition, the
penalty terms $-(\J^{(11)})(\H^{(1)})^{-1}$ and $(\J^{(22)})(\H
^{(2)})^{-1}$ are the two main diagonal blocks in the partitioned
matrix~$\Bmat$, respectively. Therefore,
\[
\tr\bigl[\bigl(\J^{(22)}\bigr) \bigl(\H^{(2)}
\bigr)^{-1}\bigr]-\tr\bigl[\bigl(\J^{(11)}\bigr) \bigl(\H
^{(1)}\bigr)^{-1}\bigr]=\tr\Bmat.
\]
We claim that the number of non-zero eigenvalues $\lambda_i$ of $\Bmat$
is $m$. To verify this, the characteristic equation $|\Bmat-\lambda
\Imat_{p_1+p_2}|=0$ can be written as
\[
\left|\matrix{ \J_{\theta\theta}+\lambda
\H_{\theta\theta}& \J_{\theta\theta} & \J _{\theta\zeta}
\vspace*{2pt}\cr
\J_{\theta\theta}& \J_{\theta\theta}-\lambda\H_{\theta\theta} & \J
_{\theta\zeta}-\lambda\H_{\theta\zeta}
\vspace*{2pt}\cr
\J_{\zeta\theta}& \J_{\zeta\theta}-\lambda\H_{\theta\zeta} &
\J_{\zeta
\zeta}-\lambda\H_{\zeta\zeta}}\right| =0.
\]
Subtract the second column from the first column, and then
subtract the first row from the second row,
\[
0=\left|
\matrix{ \lambda\H_{\theta\theta} &
\J_{\theta\theta}& \J_{\theta\zeta
}
\vspace*{2pt}\cr
\zero& -\lambda\H_{\theta\theta} & -\lambda\H_{\theta\zeta}
\vspace*{2pt}\cr
\lambda\H_{\zeta\theta} & \J_{\zeta\theta}-\lambda\H_{\zeta
\theta}&
\J_{\zeta\zeta}-\lambda\H_{\zeta\zeta} }\right|
=(-1)^{p_1}\lambda^{2p_1}\left|\matrix{ \H_{\theta\theta} & \J_{\theta\theta}& \J_{\theta\zeta}
\vspace*{2pt}\cr
\zero& \H_{\theta\theta} & \H_{\theta\zeta}
\vspace*{2pt}\cr
\H_{\zeta\theta} & \J_{\zeta\theta}-\lambda\H_{\zeta\theta}&
\J_{\zeta
\zeta}-\lambda\H_{\zeta\zeta} }\right|.
\]
If AIC is considered, the $\J$ matrices are the same as the $\H$ matrices.
Subtract the second row from the first, and then
subtract the second column multiplied by
$\H_{\theta\theta}^{-1}\H_{\theta\zeta}$ from the third to get:
\[
0 %=\lambda^{2p_1}
=\lambda^{2p_1} \left|
\matrix{\H_{\theta\theta}& \zero& \zero
\vspace*{2pt}\cr
\zero& \H_{\theta\theta}& \zero
\vspace*{2pt}\cr
\H_{\zeta\theta} & (1-\lambda)\H_{\zeta\theta}& (1-\lambda) \bigl(\H
_{\zeta
\zeta}-\H_{\zeta\theta}\H^{-1}_{\theta\theta}
\H_{\theta\zeta}\bigr)} \right|.
\]
The eigenvalues are $\lambda=0$ (multiplicity $=2p_1$) and $\lambda=1$
(multiplicity $=m$).

\textit{Proof of} (2). The required result is a direct
consequence of (1) and Lemma \ref{prob}.

\textit{Proof of} (3). For CLAIC, we show that $P^{\mathrm{CLAIC}}_1$ is
asymptotically equivalent to a non-central chi-square probability. Note
that CLAIC selects model 1 if
the CLAIC comparison is:
%
%e7.3 #&#
%
\begin{equation}
\Pr \bigl[2 \bigl\{L^{(2)}_{\mathrm{CL}}(\ghatbf_{n})-
L^{(1)}_{\mathrm{CL}}(\thatbf_{1n}) \bigr\} < 2 \bigl\{\tr
\bigl(\J^{(2)} \bigl[\H^{(2)}\bigr]^{-1}\bigr) - \tr
\bigl(\J^{(1)} \bigl[\H^{(1)}\bigr]^{-1}\bigr) \bigr\}
\bigr]. \label{eq-claicpr}
\end{equation}
Here $2[L^{(2)}_{\mathrm{CL}}(\ghatbf_{n})-L^{(1)}_{\mathrm{CL}}(\thatbf_{1n})]$ is a
non-negative quadratic form, and a representation for it is obtained below.

Write $L_{\mathrm{CL}}(\thatbf_{2n},\zehatbf_n)=L^{(2)}_{\mathrm{CL}}(\ghatbf_{n})$
and $L_{\mathrm{CL}}(\thtilbf_n(\zebf^*),\zebf^*)=L^{(1)}_{\mathrm{CL}}(\thatbf
_{1n})$, where
$\zebf^*=\zero$.
Let $\thtilbf_n(\zebf)$ be the maximum composite likelihood estimate when
$\zebf$ is fixed, so that
$L_{\mathrm{CL}}(\thtilbf_n(\zebf),\zebf)$ is the profile composite log-likelihood.

Assume that all of the regularity conditions for maximum likelihood apply
to all of the marginal densities in the composite likelihood.
The derivation below is similar to a result in Cox and Hinkley (\cite{CH}, Section~9.3) for the full log-likelihood.
For the difference of composite log-likelihoods in (\ref{eq-claicpr}),
we take
an expansion to second order:
%
%e7.4 #&#
%
\begin{eqnarray}\label{eq-D1}
&&2\bigl[L_{\mathrm{CL}}(\thatbf_{2n},\zehatbf_n)-L_{\mathrm{CL}}
\bigl(\thtilbf_n\bigl(\zebf ^*\bigr),\zebf^*\bigr)\bigr]\nonumber\\
&&\quad= n\bigl(
\thatbf_{2n}-\thbf^*\bigr)^T\H_{\th\th}\bigl(
\thatbf_{2n}-\thbf^*\bigr)
\nonumber
\\[-8pt]
\\[-8pt]
\nonumber
&&\qquad{}+2n\bigl(\thatbf_{2n}-\thbf^*\bigr)^T\H_{\th\ze}
\bigl(\zehatbf_n-\zebf^*\bigr) +n\bigl(\zehatbf_n-\zebf^*
\bigr)^T\H_{\ze\ze}\bigl(\zehatbf_n-\zebf^*\bigr)
\\
&&\qquad{}-n \bigl(\thtilbf_n\bigl(\zebf^*\bigr)-\thbf^*\bigr)^T
\H_{\th\th}\bigl(\thtilbf_n\bigl(\zebf ^*\bigr)-\thbf^*\bigr)
+\mathrm{o}_p(1).\nonumber
\end{eqnarray}
For the profile likelihood,
by differentiating ${\p L_{\mathrm{CL}}(\thtilbf_n(\zebf),\zebf)/\p\thbf
}=\zero
$, one gets:
\[
{\p^2 L_{\mathrm{CL}} \over\p\thbf\,\p\thbf^T}\bigl(\thtilbf_n(\zebf),\zebf\bigr)
{\p
\thtilbf_n\over
\p\zebf^T} + {\p^2 L_{\mathrm{CL}} \over\p\thbf\,\p\zebf^T}\bigl(\thtilbf _n(
\zebf ),\zebf\bigr)=\zero,
\]
so that as $n\to\oo$,
\[
{\p\thtilbf_n(\zebf)\over\p\zebf^T}\Bigm|_{\zebf^*} = -\H_{\th\th}^{-1}
\H_{\th\ze} +\mathrm{o}_p(1).
\]
Expand $\thtilbf_n(\zebf)$ around $\zebf=\zebf^*$ at $\zebf
=\zehatbf_n$
to get
\[
\thtilbf_n\bigl(\zebf^*\bigr)=\thtilbf_n(
\zehatbf_n)+ {\p\thtilbf_n(\zebf)\over\p\zebf^T}\Bigm|_{\zebf^*} \bigl(
\zebf^*-\zehatbf_n\bigr)+\mathrm{o}_p(1) =\thatbf_{2n} +
\H_{\th\th}^{-1}\H_{\th\ze}\bigl(\zehatbf_n-
\zebf ^*\bigr) + \mathrm{o}_p(1).
\]
Hence,
\[
\thtilbf_n\bigl(\zebf^*\bigr)-\thbf^* =\thtilbf_n\bigl(
\zebf^*\bigr)-\thatbf _{2n}+\thatbf _{2n}-\thbf^* =
\H_{\th\th}^{-1}\H_{\th\ze}\bigl(\zehatbf_n-
\zebf^*\bigr) + \bigl(\thatbf _{2n}-\thbf^*\bigr) +\mathrm{o}_p(1).
\]
Substitute into (\ref{eq-D1}) to get
%
%e7.5 #&#
%
\begin{eqnarray}\label{eq-D2}
&& 2\bigl[L_{\mathrm{CL}}(\thatbf_{2n},\zehatbf_n)-L_{\mathrm{CL}}
\bigl(\thtilbf_n\bigl(\zebf^*\bigr),\zebf ^*\bigr)\bigr]
\nonumber
\\[-8pt]
\\[-8pt]
\nonumber
&&\quad= n\bigl(
\zehatbf_n-\zebf^*\bigr)^T \bigl[\H_{\ze\ze}-
\H_{\ze\th}\H_{\th\th
}^{-1}\H_{\th
\ze}\bigr] \bigl(
\zehatbf_n-\zebf^*\bigr) + \mathrm{o}_p(1).
\end{eqnarray}
Under a sequence of contiguous alternatives,
$n^{1/2}(\E[\zehatbf_n]-\zebf_n)\to\zero$ and
$n^{1/2}(\zebf_n-\zebf^*)\to\epsbf$ as $n\to\oo$.
So marginally $n^{1/2}(\zehatbf_n-\zebf^*)$ is asymptotically
$N( \dbf_\ze,\V_\ze)$, where $ \dbf_\ze=\epsbf$
and $\V_\ze$ is the $(2,2)$ block of the partitioned covariance matrix,
\[
\pmatrix{\H^{\th\th} & \H^{\th\ze}
\cr
\H^{\ze\th}
& \H^{\ze
\ze
} } %
^{-1} %
\pmatrix{
\J_{\th\th} & \J_{\th\ze}
\cr
\J_{\ze\th} & \J_{\ze
\ze
} }
\pmatrix{\H^{\th\th} & \H^{\th\ze}
\cr
\H^{\ze\th} &
\H^{\ze
\ze
} } %
^{-1}.
\]
Then, (\ref{eq-D2}) is asymptotically a quadratic form based on
a random vector with $N( \dbf_\ze,\V_\ze)$ distribution.

For CLBIC, the arguments are similar to that of CLAIC. Here, we
highlight the differences between CLBIC and CLAIC. The result is
established based on the following comparison
%
%e7.6 #&#
%
\begin{equation}
2 \bigl\{L^{(2)}_{\mathrm{CL}}(\ghatbf_{n})-
L^{(1)}_{\mathrm{CL}}(\thatbf_{1n}) \bigr\} < %
\cases{ \log n \bigl\{\tr\bigl(\J^{(2)} \bigl[\H^{(2)}
\bigr]^{-1}\bigr) - \tr\bigl(\J^{(1)} \bigl[\H^{(1)}
\bigr]^{-1}\bigr) \bigr\}, & \quad$\mathrm{\mathrm{CLBIC}}$, \vspace *{2pt}
\cr
2 \bigl\{
\tr\bigl(\J^{(2)} \bigl[\H^{(2)}\bigr]^{-1}\bigr) -
\tr\bigl(\J^{(1)} \bigl[\H^{(1)}\bigr]^{-1}\bigr)
\bigr\}, & \quad $\mathrm{\mathrm{CLAIC}}.$} %
\end{equation}
The left-hand side has order $\mathrm{O}_p(1)$. For CLAIC, the right-hand side
is just $\mathrm{O}_p(1)$, so there is positive probability that CLAIC selects
model 2. On the contrary, for CLBIC, the right-hand side is $\mathrm{O}_p(\log
n)$. Together with the asymptotic positiveness of the penalty term
difference (see Lemma \ref{mono}), the increase in the likelihood is
offset by the increase in the penalty. Therefore, asymptotically CLBIC
cannot select model 2.\vadjust{\goodbreak}

\textit{Proof of} (4). It is similar to the proof of (3) and is
omitted here.

%s7.4 #&#
\subsection{Proof of Theorem \texorpdfstring{\protect\ref{main.robust}}{3.2}}\label{appa4}

This is similar to the proof of Theorem \ref{main}.

%s7.5 #&#
\subsection{Technical lemmas}\label{appa5}

%le7.1 #&#
\begin{lemma}\label{prob} (1) Let $Z^2_1,Z^2_2,\ldots,Z^2_m$ be independent $\chi
^2_1$ random variables. Suppose that $m'<m$. Then,
\[
P\bigl(Z^2_1+\cdots+Z^2_{m'}<2m'
\bigr)<P\bigl(Z^2_1+\cdots+Z^2_m<2m
\bigr).
\]
(2) Further let $\lambda_1,\lambda_2,\ldots,\lambda_m$ be non-negative
constants. Then,
\[
P\bigl(\lambda_1 Z^2_1+\cdots+
\lambda_m Z^2_m<2(\lambda_1+
\cdots +\lambda _m)\bigr)\leq P\bigl(Z^2_1+
\cdots+Z^2_m<2m\bigr).
\]
The equality sign holds if and only if $\lambda_1=\lambda_2=\cdots
=\lambda_m$.
\end{lemma}

\begin{pf} (1) Let $\bar{U}_m$ be the sample average of
$Z^2_1,\ldots,Z^2_m$. Below, we compare the probabilities $P(\bar
{U}_m>2)$ and $P(\bar{U}_{m+1}>2)$. It can be checked that
\[
P(\bar{U}_m>2)=\int^{\infty}_2g_m(t)\,\mathrm{d}t=
\frac{m}{2^{m/2}\Gamma
(m/2)}\int^{\infty}_2(mt)^{m/2-1}\mathrm{e}^{-mt/2}\,\mathrm{d}t
.
\]
Consider the ratio between the integrands $g_m(t)$ and $g_{m+1}(t)$,
\[
R(t)=\frac{g_{m+1}(t)}{g_m(t)}=\frac{(m+1)^{(m+1)/2}\Gamma
(m/2)}{\sqrt {2}m^{m/2}\Gamma((m+1)/2)}\sqrt{t}\mathrm{e}^{-t/2}.
\]
Note that $\sqrt{t}\mathrm{e}^{-t/2}$ is monotonic decreasing for $t>2$, it
suffices to show that $R(2)<1$. To achieve that, the Binet's formula
(see Sasv\'ari \cite{Sa}) can be employed,
\[
\frac{\Gamma(m/2)}{\Gamma((m+1)/2)}=\sqrt{2}\mathrm{e}^{1/2}\frac
{(m-2)^{(m-1)/2}}{(m-1)^{m/2}}\exp\bigl[\theta
\bigl((m-2)/2\bigr)-\theta\bigl((m-1)/2\bigr)\bigr],
\]
where
\[
\theta(x)=\int^{\infty}_0 \biggl(
\frac{1}{\mathrm{e}^t-1}-\frac{1}{t}+\frac{1}{2} \biggr)\mathrm{e}^{-xt}
\frac
{1}{t}\,\mathrm{d}t.
\]
The following bound is also used (see Lemma 2 of Sasv\'ari \cite{Sa}); for
$x>0$,
\[
\theta(x)-\theta(x+1/2)<\theta(x)-\theta(x+1)= \biggl(x+\frac
{1}{2}
\biggr)\log \biggl(1+\frac{1}{x} \biggr)-1.
\]
Then
\begin{eqnarray*}
R(2)&\leq& \mathrm{e}^{-3/2}\sqrt{2}\frac{m}{m-1} \biggl( 1+
\frac{1}{m} \biggr)^{m/2} \biggl( 1+\frac{1}{m-2}
\biggr)^{-(m-2)/2} \biggl( 1+\frac{2}{m-2} \biggr)^{(m-2)/2}
\\
&\le&\sqrt{2} \frac{m^{1/2}(m+1)^{1/2}}{m-1} \biggl( 1+\frac{1}{m-2}
\biggr)^{-(m-2)/2}.
\end{eqnarray*}
The right-hand side is monotonic decreasing series of $m$ converging to
$\sqrt{2}\mathrm{e}^{-1/2}\approx0.8578<1$. It is smaller than $1$ when
$m\geq12$. We complete the proof by reporting the numerical values
of $P(\bar{U}_m>2)=P(Z^2_1+\cdots+Z^2_m>2m)$ for $m=1,\ldots,12$. One
can see that the monotonic decreasing pattern also holds for $m\leq12$.\vspace*{9pt}

\begin{center}
\begin{tabular}{@{}lllllll@{}}\hline
$m$ &1&2&3&4&5&6\\\hline
$P(\bar{U}_m>2)$&0.157&0.135&0.112&0.092&0.075&0.062\\\hline
$m$ &7&8&9&10&11&12\\\hline
$P(\bar{U}_m>2)$&0.051&0.042&0.035&0.029&0.024&0.020 \\\hline
\end{tabular}
\end{center}
\vspace*{9pt}

(2) Let $\Omega$ be the event
$\{\lambda_1 Z^2_1+\cdots+\lambda_m Z^2_m<c(\lambda_1+\cdots
+\lambda
_m)\}$
and
\[
G(\lmbf)%=G(\lambda_1,\lambda_2,\ldots,\lambda_m)
=P\bigl(\lambda_1 Z^2_1+
\cdots+\lambda_m Z^2_m<c(\lambda_1+
\cdots +\lambda _m)\bigr)=P(\Omega),
\]
where the constant $c$ is 2. Without loss of generality,
fix the value of $\lambda_1+\cdots+\lambda_m=m$, and let
$G^*(\lambda_2,\ldots,\lambda_m)=G(m-\lambda_2-\cdots-\lambda
_m,\lambda
_2,\ldots,\lambda_m)$,
which we abbreviate below as $G^*(\lmbf)$.
We will consider (i) the stationary points of $G^*$ and (ii) boundary
points of $G^*$.

 First, we give the first-order conditions for the stationary
points. Rewrite
\begin{eqnarray*}
G^*(\lmbf)&=&K\int^{cm}_0\int
^{cm-v_m}_0\cdots\int^{cm-v_m-\cdots
-v_2}_0
\Biggl\{\prod^m_{k=1}\lambda
^{-1/2}_kv^{-1/2}_k\mathrm{e}^{-v_k/2\lambda
_k}
\Biggr\}\,\mathrm{d}v_1\cdots \,\mathrm{d}v_m.
\end{eqnarray*}
Here, $K$ is a proportionality constant. Let
\begin{eqnarray*}
E_i(\lmbf)&=&K\int^{cm}_0\int
^{cm-v_m}_0\cdots\int^{cm-v_m-\cdots
-v_2}_0v_i
\Biggl\{\prod^m_{k=1}\lambda
^{-1/2}_kv^{-1/2}_k\mathrm{e}^{-v_k/2\lambda
_k}
\Biggr\}\,\mathrm{d}v_1\cdots \,\mathrm{d}v_m,
\\
E_{ij}(\lmbf)&=&K\int^{cm}_0\int
^{cm-v_m}_0\cdots\int^{cm-v_m-\cdots
-v_2}_0v_iv_j
\Biggl\{\prod^m_{k=1}\lambda
^{-1/2}_kv^{-1/2}_k\mathrm{e}^{-v_k/2\lambda_k}
\Biggr\}\,\mathrm{d}v_1\cdots \,\mathrm{d}v_m.
\end{eqnarray*}
Differentiating $G^*(\lmbf)$ with respect to $\lambda_i$ for $i\neq
1
$, we have
\[
-\frac{1}{2\lambda_i}G^*(\lmbf)+\frac{1}{2\lambda^2_i}E_{i}(\lmbf )=-
\frac{1}{2\lambda_1}G^*(\lmbf)+\frac{1}{2\lambda
^2_1}E_{1}(\lmbf)=\nu, \qquad i=2,
\ldots,m,
\]
where $\nu$ is the Lagrange multiplier. To simplify the first-order
conditions, it is convenient to introduce the following notation. Define
\begin{eqnarray*}
h_1(\lmbf)&=&K\int^{cm}_0\int
^{cm-v_m}_0\cdots\int^{cm-v_m-\cdots
-v_3}_0
\lambda^{-3/2}_1(cm-v_m-\cdots-v_2)^{1/2}
\\
&&\hspace*{139pt}{}\times \mathrm{e}^{-(cm-v_m-\cdots-v_2)/2\lambda_1} \\
&&\hspace*{139pt}{}\times\Biggl\{\prod^m_{k=2}
\lambda ^{-1/2}_kv^{-1/2}_k\mathrm{e}^{-v_k/2\lambda_k}
\Biggr\}\,\mathrm{d}v_2\cdots \,\mathrm{d}v_m,
\\
h_{11}(\lmbf)&=&K\int^{cm}_0\int
^{cm-v_m}_0\cdots\int^{cm-v_m-\cdots
-v_3}_0
\lambda^{-5/2}_1(cm-v_m-\cdots-v_2)^{3/2}
\\
&&\hspace*{139pt}{}\times \mathrm{e}^{-(cm-v_m-\cdots-v_2)/2\lambda_1}\\
&&\hspace*{139pt}{}\times \Biggl\{\prod^m_{k=2}
\lambda ^{-1/2}_kv^{-1/2}_k\mathrm{e}^{-v_k/2\lambda_k}
\Biggr\}\,\mathrm{d}v_2\cdots \,\mathrm{d}v_m,
\\
h_{12}(\lmbf)&=&K\int^{cm}_0\int
^{cm-v_m}_0\cdots\int^{cm-v_m-\cdots
-v_3}_0
\lambda^{-3/2}_1\lambda^{-3/2}_2(cm-v_m-
\cdots -v_2)^{1/2}v^{1/2}_2
\\
&&\hspace*{139pt}{}\times \mathrm{e}^{-(cm-v_m-\cdots-v_2)/2\lambda_1}\mathrm{e}^{-v_2/2\lambda_2}\\
&&\hspace*{139pt}{}\times \Biggl\{\prod
^m_{k=3}\lambda^{-1/2}_kv^{-1/2}_k\mathrm{e}^{-v_k/2\lambda
_k}
\Biggr\}\,\mathrm{d}v_2\cdots \,\mathrm{d}v_m.
\end{eqnarray*}
Similarly, define $h_i$, $h_{ii}$, and $h_{ij}$ for other $i$,
$j$. Below are some useful results obtained from integration by parts
over variable $v_1$,\vspace*{-1pt}
\begin{eqnarray*}
E_1(\lmbf)&=&\lambda_1G^*(\lmbf)-2
\lambda^2_1h_1(\lmbf),
\\[-1pt]
E_2(\lmbf)&=&\lambda_2G^*(\lmbf)-2\lambda^2_2h_2(
\lmbf),
\\[-1pt]
E_{11}(\lmbf)&=&3\lambda_1E_1(\lmbf)-2
\lambda^3_1h_{11}(\lmbf),
\\[-1pt]
E_{12}(\lmbf)&=&\lambda_1E_2(\lmbf)-2
\lambda^2_1\lambda _2h_{12}(\lmbf)
.
\end{eqnarray*}
Then, the first order conditions becomes $h_1=h_2=\cdots=h_m=-\nu$.

Next, we show that stationary points of $G^*(\lmbf)$ without
satisfying $\lambda_1=\cdots=\lambda_m$ do not have semi-negative
definite Hessian matrix. Differentiating $G^*(\lmbf)$ with respect to
$\lambda_i$ twice,\vspace*{-1pt}
\begin{eqnarray*}
\frac{\partial^2 G^*(\lmbf)}{\partial\lambda^2_i} &=&\frac{1}{4\lambda^2_i} \biggl( 3G-\frac{6 E_i}{\lambda_i}+
\frac
{E_{ii}}{\lambda^2_i} \biggr) +\frac{1}{4\lambda^2_1} \biggl( 3G-\frac{6 E_1}{\lambda_1}+
\frac
{E_{11}}{\lambda^2_1} \biggr)
\\[-1pt]
&&{}+\frac{2}{4\lambda_i\lambda_1} \biggl( G-\frac{E_i}{\lambda
_i}-\frac
{E_1}{\lambda_1}+
\frac{E_{1i}}{\lambda_i\lambda_1} \biggr)
\\[-1pt]
&=&\frac{h_i}{\lambda_i}+\frac{h_1}{\lambda_1}+\frac{1}{2\lambda
_i}(h_i-h_{ii})+
\frac{1}{2\lambda_1}(h_1-h_{11})-\frac{1}{\lambda
_i}(h_1-h_{1i})
.
\end{eqnarray*}
Below, we see that the right-hand side must be positive if $\lambda
_1\neq\lambda_i$ and therefore cannot be a local maximum. By
definitions, the first two terms are positive. For the third term,
consider the quantities defined below,
\[
R_{ii}=G^*(\lmbf)-\frac{2E_i(\lmbf)}{\lambda_i}+\frac{E_{ii}(\lmbf
)}{\lambda^2_i}.
\]
It can be rewritten as the integration of the product of
$(1-v_i/\lambda
_i)^2$ and some positive terms. Therefore, $R_{ii}$ must be positive.
In addition, $R_{ii}=2\lambda_i(h_i-h_{ii})$. Then, we show that
$h_i-h_{ii}>0$. The fourth term can be handled in the same manner.
For the last term, the symmetry $E_{1i}=E_{i1}$ implies
$\lambda_i h_i+\lambda_1 h_{1i} = \lambda_1 h_1+\lambda_i h_{i1}$; then
using the first order condition for a stationary point and the symmetry
of $h_{ij}$, $(\lambda_1-\lambda_i)(h_1-h_{1i})=0$.
If $\lambda_1\neq\lambda_i$, then $h_1-h_{1i}=0$. The stationary
point must not be a local maximum.

 Now, we have shown that $\lmbf=\one_m$
is the only stationary point of $G^*(\lmbf)$ that could be a local
maximum. It should be noted that such stationary point is not
necessarily a local maximum. To avoid the difficulties in checking the
negative-definiteness of the Hessian matrix, an indirect approach is
adopted. Here, we compare the unique stationary point with the boundary
points. The boundary is defined by $\{\lambda_i=0\mbox{ for some
}i=1,2,\ldots,m\}$. Result (2) on the boundary points can be
established by applying result (1) and result (2) for stationary points
inductively. (Note: for any $c$, $\lmbf=\one_m$ is always a
stationary point. However, result (1) is not necessarily valid for all~$c$,
so, the local maximality does not always hold for any~$c.$)
\end{pf}
%
%le7.2 #&#
\begin{lemma}[(Monotonicity of the penalty term $\tr(\H^{-1}\J)$)]\label{mono}
If model 1 is nested within model 2, $\tr((\H^{(1)})^{-1}\J
^{(1)})<\tr
((\H^{(2)})^{-1}\J^{(2)})$, if $\H^{(1)},\J^{(1)}$ are evaluated at
$\thbf^*$ and
$\H^{(2)},\J^{(2)}$ are evaluated at $(\thbf^*,\zebf^*)$.
\end{lemma}

\begin{pf}
Suppose that the parameters are $(\zebf
^*,\thbf
^*)$ for model 2
and $\thbf^*$ for model 1. Below, if not specified, the arguments
of the $\H,\J$ matrices are $(\zebf^*,\thbf^*)$.
For model 1, the penalty term is $\tr(\H^{-1}_{\theta\theta}\J
_{\theta
\theta})$.

Next, we consider the partitioning:\vspace*{1pt}
\[
\left( %
\matrix{ \H_{\zeta\zeta}&\H_{\zeta\theta}
\vspace*{2pt}\cr
\H_{\theta\zeta}&\H_{\theta\theta} }
\right)^{-1} \left( \matrix{ \J_{\zeta\zeta}&
\J_{\zeta\theta}
\vspace*{2pt}\cr
\J_{\theta\zeta}&\J_{\theta\theta} }
\right).
\]
We have (see Morrison \cite{Mo}, Section~2.11)\vspace*{3pt}
\begin{eqnarray*}
\left( %
\matrix{ \H_{\zeta\zeta}&
\H_{\zeta\theta}
\vspace*{2pt}\cr
\H_{\theta\zeta}&\H_{\theta\theta} }
\right)^{-1}&=& \left( %
\matrix{ \bigl(\H_{\zeta\zeta}-
\H_{\zeta\theta}\H^{-1}_{\theta\theta}\H _{\theta\zeta}
\bigr)^{-1}\vspace*{2pt}\cr
-\H^{-1}_{\theta\theta}\H_{\theta\zeta}\bigl(\H_{\zeta\zeta}-\H
_{\zeta\theta
}\H^{-1}_{\theta\theta}\H_{\theta\zeta}
\bigr)^{-1}}\right.\\[3pt]
&&\hspace*{7pt}\left.\matrix{
 -\bigl(\H_{\zeta\zeta}-\H_{\zeta\theta}
\H^{-1}_{\theta\theta}\H _{\theta\zeta
}\bigr)^{-1}
\H_{\zeta\theta}\H^{-1}_{\theta\theta}
\vspace*{2pt}\cr
 \H^{-1}_{\theta\theta}+\H^{-1}_{\theta\theta}
\H_{\theta\zeta
}\bigl(\H_{\zeta
\zeta}-\H_{\zeta\theta}\H^{-1}_{\theta\theta}
\H_{\theta\zeta
}\bigr)^{-1}\H _{\zeta\theta}\H^{-1}_{\theta\theta}}
\right).
\end{eqnarray*}
The change in the penalty term is therefore
\begin{eqnarray*}
&&\tr\bigl(\H_{\zeta\zeta}-\H_{\zeta\theta}\H^{-1}_{\theta
\theta}
\H _{\theta\zeta}\bigr)^{-1} \bigl(\J_{\zeta\zeta}-\H_{\zeta\theta}
\H^{-1}_{\theta\theta}\J _{\theta\zeta
}\bigr)
\\
&&\quad{}+\tr\H^{-1}_{\theta\theta}\H_{\theta\zeta}\bigl(\H_{\zeta\zeta
}-
\H_{\zeta
\theta}\H^{-1}_{\theta\theta}\H_{\theta\zeta}
\bigr)^{-1}\bigl(\H_{\zeta
\theta}\H ^{-1}_{\theta\theta}
\J_{\theta\theta}-\J_{\zeta\theta}\bigr)
\\
&&\qquad=\tr\bigl(\H_{\zeta\zeta}-\H_{\zeta\theta}\H^{-1}_{\theta\theta
}
\H_{\theta
\zeta}\bigr)^{-1} \bigl(\J_{\zeta\zeta}-\H_{\zeta\theta}
\H^{-1}_{\theta\theta}\J _{\theta\zeta} -\J_{\zeta\theta}
\H^{-1}_{\theta\theta} \H_{\theta\zeta} +\H_{\zeta\theta}
\H^{-1}_{\theta\theta}\J_{\theta\theta}\H ^{-1}_{\theta
\theta}
\H_{\theta\zeta}\bigr).
\end{eqnarray*}
Note that the term
\[
\J_{\zeta\zeta}-\H_{\zeta\theta}\H^{-1}_{\theta\theta}\J
_{\theta\zeta} -\J_{\zeta\theta} \H^{-1}_{\theta\theta}
\H_{\theta\zeta} +\H_{\zeta\theta}\H^{-1}_{\theta\theta}
\J_{\theta\theta}\H ^{-1}_{\theta
\theta}\H_{\theta\zeta}
\]
must be positive definite because $\J$ has the form $\E[\nabla\nabla
^T]$. It is the same as
\[
\E\bigl[\nabla_{\zeta}-\H_{\zeta\theta}\H^{-1}_{\theta\theta}
\nabla _{\theta
}\bigr] \bigl[\nabla_{\zeta}-\H_{\zeta\theta}
\H^{-1}_{\theta\theta}\nabla _{\theta
}\bigr]^T,
\]
where $\nabla_{\theta}$ and $\nabla_{\zeta}$ are the gradients of the
composite log-likelihood with respect to $\thbf$ and $\zebf$, respectively.
The term $(\H_{\zeta\zeta}-\H_{\zeta\theta}\H^{-1}_{\theta
\theta}\H
_{\theta\zeta})^{-1}$ is also positive definite because it is a
principal block of the matrix
\[
\left( %
\matrix{ \H_{\zeta\zeta}&\H_{\zeta\theta}
\vspace*{2pt}\cr
\H_{\theta\zeta}&\H_{\theta\theta} }
\right)^{-1}.
\]
For any two positive definite matrices $A$ and $B$, the trace $\tr
(AB)$ must be positive. To see this, consider eigenvalue decomposition
$A=P\Lambda P^T$. The trace $\tr(AB)=\tr(\Lambda P^TBP)$ is the dot
product of the diagonals of $\Lambda$ and $P^TBP$. Since $P^TBP$ is
positive definite, all diagonal elements must be positive. We have the
desired results that the penalty term is monotonic increasing.
\end{pf}
%s8 #&#
\section{Full and composition likelihoods of the linear mixed-effects model}\label{appb}

For the multivariate normal mixed-effects model Laird and Ware \cite
{LW}, both
the full likelihood
and composite likelihood can be computed readily, after making use of
results on vec and vech operations (see Fackler \cite{Fa}, Magnus and
Neudecker \cite{MN}).

\textit{Model}:
\begin{eqnarray*}
\Ybf_i&=&\xbf_i\bebf+\mathbf{z}_i
\bbf_i+\epsbf_i,\qquad i=1,2,\ldots,n,
\\
\bbf_i&\sim &N(\zero,\Psi), \qquad\epsbf_i\sim N(\zero,\phi
\Imat _d),
\end{eqnarray*}
where $\bebf$ is $(s+1)$-dimensional vector of fixed effects, $\bbf
_i$ is
$r$-dimensional vector of random effects. $\xbf_i$ and $\mathbf{z}_i$ are
$d\times s$ and $d\times r$ observable matrices, $\xbf_i$
has a first column of 1s,
$\phi$ is a variance parameter, $\Psi$ is a $r\times r$ covariance matrix.

\textit{Conventions}: Define the duplication matrix $\Dmat_r$
such that for any $r\times r$ symmetric matrix~$\Amat$, we have $\vec
\Amat=\Dmat_r \vech\Amat$. Define permutation matrices $\Tmat_{rr}$
such that for any $r\times r$ matrix~$\Amat$, we have $\Tmat_{rr}
\vec\Amat= \vec\Amat^T$. Define the duplication matrix $\Dmat_r$
and elimination matrix $\Emat_r$ such that for any $r\times r$
symmetric matrix $\Amat$, we have $\Emat_r \vec\Amat=\vech\Amat$
and $\vec\Amat=\Dmat_r \vech\Amat$.
The duplication matrix is unique but not the elimination matrix;
for the latter, it is convenient to operate on the lower triangle.
Let $\Imat_r$ be the $r\times r$ identity matrix.
Some properties of the above-mentioned
matrices are as follows. (1) $(\Imat_r+\Tmat_{rr})\Dmat_r=2\Dmat_r$,
$\Dmat_r\Emat_r(\Imat_r+\Tmat_{rr})=\Imat_r+\Tmat_{rr}$. (2) If
$\Cmat
$ is lower-triangular, we have $\vec\Cmat=\Emat^T_r \vech\Cmat$.

Details for the full likelihood and the pairwise composite likelihood are
given in two subsections below. The ideas are similar for other
composite likelihoods.

%s8.1 #&#
\subsection{Full likelihood}\label{appb1}

Define
\[
\Omega_i=\mathbf{z}_i\Psi\mathbf{z}^\Tmat_i+
\phi\Imat\quad\mbox {and}\quad \Smat_i=(\mathbf{y}_i-
\xbf_i\bebf) (\mathbf{y}_i-\xbf_i
\bebf)^T.
\]
The likelihood function is
\[
L(\beta,\Psi,\phi)=\sum^n_{i=1}
\ell_i(\bebf,\Psi,\phi;\mathbf{ y}_i,\xbf_i),
\]
where
\[
\ell_i(\beta,\Psi,\phi)=\ell_i(\bebf,\Psi,\phi;\mathbf
{y}_i,\xbf_i)= -\tfrac{1}{2} \bigl\{ \tr\bigl(
\Omega^{-1}_i \Smat_i\bigr)+\log|
\Omega_i| \bigr\} - \tfrac{1}{2}\log(2\uppi).
\]
The following alternative parameterization is beneficial to numerical
computation.
Consider $\Psi=\Cmat\Cmat^T$ and $\phi=\kappa^2$, were $\Cmat$ is lower
triangular matrix.
We have
\[
\frac{\mathrm{d} \vech\Psi}{\mathrm{d} \vech\Cmat}=\Emat_r (\Imat_r+\Tmat
_{rr}) (\Cmat \otimes\Imat_r)\Emat^T_r
.
\]

Under the $(\Cmat,\kappa)$ parametrization, the score function and
Fisher's information matrix are given as follows.

\textit{Score function}:
\begin{eqnarray*}
\frac{\mathrm{d}\ell_i}{\mathrm{d}\bebf}&=&(\mathbf{y}_i-\xbf_i
\bebf)^T\Omega ^{-1}_i\xbf_i,
\\
\frac{\mathrm{d}\ell_i}{\mathrm{d} \vech C}&=&\vec^T\bigl[\mathbf{z}^T_i
\Omega ^{-1}_i(\Smat _i-\Omega
_i)\Omega^{-1}_i \mathbf{z}_i
\Cmat\bigr]E^T_{r},
\\
\frac{\mathrm{d}\ell_i}{\mathrm{d}\kappa}&=&\kappa\tr\bigl\{ \Omega^{-2}_i(\Smat
_i-\Omega _i)\bigr\}.
\end{eqnarray*}

\textit{Fisher information matrix}:
\begin{eqnarray*}
\E\frac{\mathrm{d}}{\mathrm{d}\bebf} \biggl(\frac{\mathrm{d}\ell_i}{\mathrm{d}\bebf} \biggr)&=&-\xbf
^T_i\Omega ^{-1}_i
\xbf_i,
\\
\E\frac{\mathrm{d}}{\mathrm{d} \vech\Cmat} \biggl(\frac{\mathrm{d}\ell_i}{\mathrm{d}\bebf} \biggr)&=&\zero,
\\
\E\frac{\mathrm{d}}{\mathrm{d}\kappa} \biggl(\frac{\mathrm{d}\ell_i}{\mathrm{d}\bebf} \biggr)&=&\zero,
\\
\E\frac{\mathrm{d}}{\mathrm{d} \vech\Cmat} \biggl(\frac{\mathrm{d}\ell_i}{\mathrm{d} \vech\Cmat
} \biggr)&=& -\Emat_r
\bigl(\Cmat^T\otimes\Imat_r\bigr)\bigl\{ \bigl[
\mathbf{z}^T_i\Omega ^{-1}_i
\mathbf{z} _i\bigr]\otimes\bigl[\mathbf{z}^T_i
\Omega^{-1}_i\mathbf{z}_i\bigr]\bigr\}(\Imat
_r+\Tmat _{rr}) (\Cmat \otimes\Imat_r)
\Emat^T_r,
\\
\E\frac{\mathrm{d}}{\mathrm{d}\kappa} \biggl(\frac{\mathrm{d}\ell_i}{\mathrm{d} \vech\Psi} \biggr)&=&-2\kappa
\Emat_r \bigl(\Cmat^T\otimes\Imat_r\bigr)
\vec \bigl\{ \bigl[\mathbf {z}^T_i\Omega
^{-2}_i\mathbf{z} _i\bigr] \bigr\},
\\
\E\frac{\mathrm{d}}{\mathrm{d}\kappa} \biggl(\frac{\mathrm{d}\ell_i}{\mathrm{d}\kappa} \biggr)&=&-2\kappa^2
\tr \bigl(\Omega^{-2}_i\bigr).
\end{eqnarray*}

%s8.2 #&#
\subsection{Composite likelihood}\label{appb2}

We show details of the pairwise composite log-likelihood for
the multivariate Gaussian linear mixed-effects model.
Define the composite likelihood as
\[
L_{\mathrm{CL}}(\bebf,\Psi,\phi)=\sum^n_{i=1}
\sum_{1\le j<k\le d} \log f_{jk}(y_{ij},y_{ik};
\bebf,\Psi,\phi,\xbf_i),
\]\eject\noindent
where
\begin{eqnarray*}
&&\log f_{jk}(y_{ij},y_{ik};\bebf,\Psi,\phi,
\xbf_i) =\ell_{i,jk}(\bebf,\Psi,\phi)
\\
&&\quad= -\tfrac{1}{2} \bigl\{ \tr\bigl[ \bigl({e}^T_{jk}
\Omega _ie_{jk}\bigr)^{-1}\bigl({e}^T_{jk}
\Smat _ie_{jk}\bigr)\bigr]+\log\bigl|{e}^T_{jk}
\Omega_ie_{jk}\bigr| \bigr\}-\tfrac
{1}{4}d(d-1)\log (2
\uppi).
\end{eqnarray*}
Let $\ell_{i\mathrm{CL}}=\sum_{1\le j<k\le d} \ell_{i,jk}(\bebf,\Psi,\phi)$.
For convenience, for $a=1,2,3$, define
\begin{eqnarray*}
\Amat_{ai}&=&\sum_{jk}e_{jk}
\bigl({e}^T_{jk}\Omega_ie_{jk}
\bigr)^{-a}{e}^T_{jk},
\\
\Bmat_i&=&\sum_{jk} \bigl\{
\bigl[e_{jk}\bigl({e}^T_{jk}\Omega
_ie_{jk}\bigr)^{-1}{e}^T_{jk}
\bigr]\otimes\bigl[e_{jk}\bigl({e}^T_{jk}\Omega
_ie_{jk}\bigr)^{-1}{e}^T_{jk}
\bigr]\bigr\},
\end{eqnarray*}
where $e_{jk}$ is the $d\times2$ matrix that has 1 in the
$(j,1)$ and $(k,2)$ positions and 0 elsewhere (premultiplying by ${e}^T_{jk}$
and postmultiplying by $e_{jk}$ extracts the appropriate $2\times2$
subcovariance matrix).

\textit{Score function}: With the above alternative
parameterization of $\Cmat$ and $\kappa$, we have
\begin{eqnarray*}
\frac{\mathrm{d} \ell_{i\mathrm{CL}}}{\mathrm{d}\bebf}&=&(\mathbf{y}_i-\xbf_i
\bebf)^T\Amat _{1i}\xbf_i,
\\
\frac{\mathrm{d} \ell_{i\mathrm{CL}}}{\mathrm{d }\vech\Cmat}&=&\vec^T (\Smat_i-\Omega
_i)\Bmat _i(\mathbf{z}_i\otimes
\mathbf{z}_i) (\Cmat\otimes\Imat_r)\Emat
^T_r,
\\
\frac{\mathrm{d} \ell_{i\mathrm{CL}}}{\mathrm{d}\kappa}&=&\kappa\tr\bigl\{ \Amat_{2i} (\Smat
_i-\Omega _i)\bigr\}.
\end{eqnarray*}

\textit{Second moment $\J$ of score function}:
\begin{eqnarray*}
&&\E \biggl(\frac{\mathrm{d} \ell_{i\mathrm{CL}}}{\mathrm{d}\bebf} \biggr)^T \biggl(\frac{\mathrm{d} \ell
_{i\mathrm{CL}}}{\mathrm{d}\bebf}
\biggr)= \xbf^T_i\Amat_{1i}
\Omega_i\Amat _{1i}\xbf_i,
\\
&&\E \biggl(\frac{\mathrm{d} \ell_{i\mathrm{CL}}}{\mathrm{d}\bebf} \biggr)^T \biggl(\frac{\mathrm{d} \ell
_{i\mathrm{CL}}}{\mathrm{d} \vech C}
\biggr)=\zero,
\\
&&\E \biggl(\frac{\mathrm{d} \ell_{i\mathrm{CL}}}{\mathrm{d}\bebf} \biggr)^T \biggl(\frac{\mathrm{d} \ell
_{i\mathrm{CL}}}{\mathrm{d}\kappa}
\biggr)=\zero,
\\
&&\E \biggl(\frac{\mathrm{d} \ell_{i\mathrm{CL}}}{\mathrm{d} \vech\Cmat} \biggr)^T \biggl(\frac
{\mathrm{d} \ell
_{i\mathrm{CL}}}{\mathrm{d} \vech\Cmat}
\biggr)\\
&&\quad=\Emat_r\bigl(\Cmat^T\otimes\Imat
_r\bigr) \bigl(\mathbf{z} ^T_i\otimes
\mathbf{z}^T_i\bigr)\Bmat_i(
\Omega_i\otimes\Omega_i)\Bmat _i(\mathbf{z}
_i\otimes\mathbf{z}_i) (\Imat_r+T_{rr})
(\Cmat\otimes\Imat_r) \Emat ^T_r,
\\
&&\E \biggl(\frac{\mathrm{d} \ell_{i\mathrm{CL}}}{\mathrm{d} \vech\Cmat} \biggr)^T \biggl(\frac
{\mathrm{d} \ell
_{i\mathrm{CL}}}{\mathrm{d}\kappa}
\biggr)=2\kappa\Emat_r\bigl(\Cmat^T\otimes\Imat
_r\bigr) \bigl(\mathbf{z} ^T_i\otimes
\mathbf{z}^T_i\bigr)\Bmat_i(
\Omega_i\otimes\Omega_i)\vec (\Amat_{2i}),
\\
&&\E \biggl(\frac{\mathrm{d} \ell_{i\mathrm{CL}}}{\mathrm{d}\kappa} \biggr) \biggl(\frac{\mathrm{d} \ell
_{i\mathrm{CL}}}{\mathrm{d}\kappa} \biggr)=2
\kappa^2 \tr(\Amat_{2i}\Omega_i\Amat
_{2i}\Omega_i).
\end{eqnarray*}

\textit{Expectation of Hessian matrix $\H$}:
\begin{eqnarray*}
\E\frac{\mathrm{d}}{\mathrm{d}\bebf} \biggl(\frac{\mathrm{d} \ell_{i\mathrm{CL}}}{\mathrm{d}\bebf} \biggr)&=&-\xbf
^T_i\Amat_{1i}\xbf_i,
\\
\E\frac{\mathrm{d}}{\mathrm{d} \vech\Cmat} \biggl(\frac{\mathrm{d} \ell_{i\mathrm{CL}}}{\mathrm{d}\bebf} \biggr)&=&\zero,
\\
\E\frac{\mathrm{d}}{\mathrm{d}\kappa} \biggl(\frac{\mathrm{d} \ell_{i\mathrm{CL}}}{\mathrm{d}\bebf} \biggr)&=&\zero,
\\
\E\frac{\mathrm{d}}{\mathrm{d }\vech\Cmat} \biggl(\frac{\mathrm{d} \ell_{i\mathrm{CL}}}{\mathrm{d} \vech\Cmat
} \biggr)&=&-\Emat_r
\bigl(\Cmat^T\otimes\Imat_r\bigr) \bigl(
\mathbf{z}^T_i\otimes\mathbf{ z}_i\bigr)\Bmat(
\mathbf{z} _i\otimes\mathbf{z}_i) (\Imat_r+
\Tmat_{rr}) (\Cmat\otimes\Imat_r) \Emat^T_r,
\\
\E\frac{\mathrm{d}}{\mathrm{d}\kappa} \biggl(\frac{\mathrm{d} \ell_{i\mathrm{CL}}}{\mathrm{d} \vech\Cmat
} \biggr)&=&-2\kappa
\Emat_r\bigl(\Cmat^T\otimes\Imat_r\bigr)
\vec\bigl(\mathbf {z}^T_i\Amat _{2i}\mathbf{z}
_i\bigr),
\\
\E\frac{\mathrm{d}}{\mathrm{d}\kappa} \biggl(\frac{\mathrm{d} \ell_{i\mathrm{CL}}}{\mathrm{d}\kappa} \biggr)&=&-2\kappa ^2
\tr(\Amat_{2i}).
\end{eqnarray*}
\end{appendix}

\section*{Acknowledgements}
This research has been supported by an NSERC Discovery grant.
We are grateful to the referees for comments that have led to an
improved presentation.
%
% imsref loaded by akundreckaite, 2013-08-13 08:47:41

%


%suskaldyti doi

\printhistory

\end{document}